\documentstyle[draft,12pt,amssymb]{amsart}
\theoremstyle{plain}   %% This is the default, anyway
\newtheorem{thm}{Theorem}[section]   % Numbered within each section
\newtheorem{cor}[thm]{Corollary}     % Numbered along with thm
\newtheorem{prop}[thm]{Proposition}  % Numbered along with thm
\newtheorem{ex}[thm]{Example}        % Numbered along with thm
\theoremstyle{definition}
\newtheorem{defn}[thm]{Definition}   % Numbered along with thm
\theoremstyle{remark}
\newtheorem{rem}[thm]{Remark}        % Numbered along with thm
\newtheorem{notation}{Notation}
\numberwithin{equation}{section}
 \let\phi=\varphi
 \let\epsilon=\varepsilon
\font\fivesym=cmsy5
\newcommand{\leftmapsto}%
{\longleftarrow\hspace{-3pt}\raise1.5pt\llap{\fivesym j}}
\newcommand{\rightmapsto}%
{\raise1.5pt\rlap{\fivesym j}\hspace{-3pt}\longrightarrow}
\newcommand{\Ass}{\operatorname{Ass}}
\newcommand{\Ann}{\operatorname{Ann}}
\newcommand{\Supp}{\operatorname{Supp}}
\newcommand{\Quot}{\operatorname{Quot}}
\begin{document}
\hbox to \hsize{\hfill}
\vspace{3cm}

\begin{center}
{\Large VARIOUS NOTIONS OF ASSOCIATED PRIME IDEALS}
\\[1cm]
{\Large Robert W. Berger}
\end{center}
\section*{Introduction}
In the theory of modules over commutative rings there are several
possibilities of defining associated prime ideals.
The usual definition of an associated prime ideal $\frak p$ for a
module $M$ is that $\frak p$ is the annihilator of an element of $M$.
In \cite{Bourbaki-Alg-Comm-4} \S 1 exercise 17 a generalization of
this notion is given. $\frak p$ is called weakly associated
(faiblement associ\'e) to $M$ if $\frak p$ is minimal in the set of
the prime ideals containing the annihilator of an element of $M$
(see Definition \ref{def-Ass-essential-first-kind}). In
this paper a further generalization of this notion will be given
(Definition \ref{def-assprim}). We use ideas of Krull
\cite{Krull-Ringe-ohne-Endlichkeit}.\\
As long as the modules are noetherian all these definitions are
equivalent. But for non noetherian $R$-modules this is no longer true,
even if the ring $R$ is noetherian.\\
In this paper we give a selfcontained introduction to the various
concepts and discuss their relation with the support and the radical
of a module.  Then we illustrate by examples the scope of the
notions. For a comprehensive introduction the theory we refer to the
now classic lecture notes \cite{Serre-Alg-loc} of Serre and to
\cite{Bourbaki-Alg-Comm-4}.
Another extensive exposition of the general theory with many examples
was given by Stefan Mittelbach in \cite{Mittelbach-Dipl}.
\\
Throughout this paper ``ring'' always denotes a commutative ring
with unit element denoted by $1$. If $M$ is an $R$-module we always assume
that $1\cdot x=x$ for all $x\in M$.\\
In the first section we recall some basic definitions and facts from
``additive ideal theory''.
\section{Primary Decomposition}
\subsection{Primary and Coprimary Modules. Primary Decomposition.}
Let $R$ be a ring, $M$ an $R$-module.
\begin{defn}
A submodule $F$ of $M$ is called ``indecomposable in $M$'' iff
from $F=F_1\cap F_2$,\quad $F_1,F_2$ submodules of $M$, it follows that
$F_1=F$ or $F_2=F$.
\end{defn}
\begin{rem}
\label{zero-indcomp}
Obviously $F$ is indecomposable in $M$ iff $(0)$ is indecomposable in
$M/F$.\\
It is well known that in a noetherian $R$-Module every submodule is
can be written as an intersection of finitely many indecomposable
submodules.
\end{rem}
\begin{defn}
An element $\xi\in R$ is called a ``zero divisor for M'',
iff there exists an element $0\ne x\in M$ with $\xi\cdot x=0$.\\
An element $\xi\in R$ is called ``nilpotent for $M$'',
iff for every $y\in M$ there exists a natural number $n=n(y)$
with $\xi^n\cdot y=0$.
(In \cite{Bourbaki-Alg-Comm-4} these elements are called
{\em presque nilpotent}.)
\end{defn}
\begin{rem}
\label{rem:zerodivisorsmultclosed}\strut
\begin{enumerate}
\item
The set of all non zero divisors for $M$ is a multiplicatively
closed subset of $R$.
\item
The set of all nilpotent elements for $M$ is an ideal of $R$.
\item
If $M\ne(0)$ every nilpotent element for $M$ is also a zero divisor
for $M$.
\end{enumerate}
\end{rem}
\begin{defn}(\cite{Serre-Alg-loc})
$M$ is called ``coprimary'', iff $M\ne(0)$ and every zero divisor for
$M$ is nilpotent for $M$.
\end{defn}
\begin{rem}
\label{submodule-coprimary-is coprimary}
 From Remark \ref{rem:zerodivisorsmultclosed} it follows that the zero
divisors for a coprimary module $M$ form an ideal $\ne R$, whose complement
is multiplicatively closed, i.e. a {\em prime ideal} $\frak p$.\\
We say that ``$M$ is $\frak p$-coprimary''.\\
Obviously every non zero submodule of a $\frak p$-coprimary module is
again $\frak p$-coprimary.
\end{rem}
\begin{ex}
\label{R-mod-p-is-coprimary}
Let $V$ be a cyclic $R$-module whose annihilator $\Ann_R(V)$ is a
prime ideal $\frak p$. Then $V$ is $\frak p$-coprimary.
\end{ex}
\begin{pf}
By definition we have $V=R\cdot x\cong R/\frak p$. So each element of
$\frak p$ is nilpotent for $V$. On the other hand every zero divisor
for $R/\frak p$ lies in $\frak p$, since $\frak p$ is a prime ideal.
\end{pf}
\begin{prop}[Noether]
\label{indcomp-coprim}
Let $M\ne(0)$ be a noetherian $R$-module, $(0)$ indecomposable in
$M$. Then $M$ is coprimary.
\end{prop}
\begin{pf}
(indirect) If $M$ was not coprimary, there would be a non nilpotent
zero divisor $\rho\in R$. Let
$M_i:=\{x\mid x\in M,\quad \rho^i\cdot x=0\}$.
$M_i$ is a submodule of $M$ and $M_i\subseteq M_{i+1}$. Since $M$ is
noetherian there exists an $n\in\Bbb N$ with $M_n=M_{n+1}$. Further
we have $M_1\ne(0)$ because $\rho$ is a zero divisor for $M$ and
$\rho^n M\ne(0)$ since $\rho$ is not nilpotent for $M$.
But $M_1\cap\rho^n M\stackrel{!}{=}(0)$ (and so $(0)$ would be
decomposable in $M$).\\
Proof: If $x\in M_1\cap\rho^n M$ then
$\rho x=0$ and $x=\rho^n y$ with an $y\in M$. $\Rightarrow
\rho^{n+1}y=0$ $\Rightarrow y\in M_{n+1}$. But $M_{n+1}=M_n$ and
therefore $x=\rho^n y=0$.
\end{pf}
If $M$ is non noetherian then Proposition \ref{indcomp-coprim} does
not hold. See Example \ref{0-indcomp-not-coprim}.
\begin{defn}
A submodule $N$ of $M$ is called ``$\frak p$-primary (or primary
for $\frak p$) in $M$'', iff
$M/N$ is $\frak p$-coprimary.
\end{defn}
Immediately from the definitions follows:
\begin{rem}\strut
\begin{enumerate}
\item $M$ is coprimary $\Longleftrightarrow$ $(0)$ is $\frak p$-primary in
$M$.
\item Let $F$ and $N$ be submodules of $M$ with $M\supseteq F\supseteq
N$. Then:\\
$F$ is $\frak p$-primary in $M$ $\Longleftrightarrow$ $F/N$ is
$\frak p$-primary in $M/N$.
\item If $M=R$ and $N=\frak q$ an ideal in $R$ then $\frak q$ is
primary in $R$ iff $\frak q$ is what is called a primary ideal.
\end{enumerate}
\end{rem}
With these notions Proposition \ref{indcomp-coprim} can be
reformulated as:
\begin{cor}
\label{indecomp-primary}
Let $N$ be a proper submodule of $M$ which is indecomposable in $M$
and $M/N$ noetherian. Then $N$ is primary in $M$.
\end{cor}
The converse of Corollary \ref{indecomp-primary} is not true but:
\begin{prop}
\label{intersection-for-same-primideal}
Let $N_1,\dots,N_r$ be $\frak p$-primary submodules of $M$.\\
Then $\bigcap\limits_{i=1}^r N_i$ is also $\frak p$-primary in $M$.
\end{prop}
\begin{pf}
1) Every element of $\frak p$ is nilpotent for
$N:=\bigcap\limits_{i=1}^r N_i$:\\
If $\rho\in\frak p$ and
$x\in M$ then for for every $i$ there is an $n_i\in\Bbb N$ with
$\rho^{n_i}\cdot x\in N_i$. Then with
$n:=\max\{n_1,\dots,n_r\}$ we have $\rho^n\cdot x\in N$.\\
2) Every zero divisor for $M/N$ is an element of $\frak p$:\\
Let $\rho\in R$ be a zero divisor for $M/N$. Then there exists an
$x\in M\setminus N$ with $\rho\cdot x\in N$. Then $x\notin N_{i_0}$
for an $i_0$ and $\rho\cdot x\in N_{i_0}$. Then $\rho$ is a
zero divisor for $M/N_{i_0}$ and therefore $\rho\in\frak p$.
\end{pf}
\begin{defn}
\label{primary-decomp}
Let $N$ be a submodule of $M$. A decomposition
$$
N=\bigcap\limits_{i=1}^r F_i,\qquad \text{where the } F_i \text{ are }
\frak p_i \text{-primary in }M
$$
is called a ``primary decomposition (or representation) of $N$ in
$M$''.\\
The $\frak p_i$ are called the ``prime ideals belonging to the
primary decomposition''\\
A primary decomposition of $N$ in $M$ is called ``reduced'' (or
``irredundant'') or a ``normal decomposition (or representation)
of $N$ in $M$'' if
the following two conditions hold:
\begin{enumerate}
\item
\label{unequal-primes}
$i\ne k\Longrightarrow\frak p_i\ne\frak p_k$.
\item
\label{superfluous}
No $F_i$ contains the intersection of the others
\end{enumerate}
\end{defn}
\begin{rem}
\label{normal-repr-exists}
 From a primary decomposition one can obtain a normal
representation:\\
Group together the $F_i$ which are primary for the same prime
ideal $\frak p$ and take their intersection. By Proposition
\ref{intersection-for-same-primideal} this is again a
$\frak p$-primary submodule of $M$. This take care of condition
(\ref{unequal-primes}) of Definition \ref{primary-decomp}.
Then omit a primary submodule that contains the intersection of
the others. Proceed until also condition (\ref{superfluous}) of
Definition \ref{primary-decomp} is satisfied.
\end{rem}
\begin{ex}
\label{coprimary-normal}
In a $\frak p$-coprimary $R$-module $M$ the submodule $(0)$ has the
normal representation $(0)=(0)$, since $(0)$ is a $\frak p$-primary
submodule of $M$.
\end{ex}
 From Remark \ref{zero-indcomp}, Corollary \ref{indecomp-primary}
together with Remark \ref{normal-repr-exists} we obtain the well
known
\begin{prop}
If $M$ is noetherian then every proper submodul $N$ of $M$
has a normal representation in $M$.
\end{prop}
\subsection{Quotient Modules and \boldmath $S$-Components.}
Let $M$ be an $R$-module, $S$ a multiplicatively closed subset of $R$
with $1\in S$,
$$
\psi:M\longrightarrow M_S
$$
the canonical homomorphism of M into the quotient module of $M$ with
respect to $S$.
\begin{notation}
For an $R_S$-submodule $U$ of $M_S$ we denote by $U\cap M$ the full inverse
image of $U$ under $\psi$:
$$
U\cap M:=\psi^{-1}(U).
$$
$U\cap M$ obviously is an $R$-submodule of $M$.\\[1ex]
For an $R$-submodule $N$ of $M$ we denote by $R_S\cdot N$ the
$R_S$-submodule of $M_S$ generated by $\psi(N)$:
$$
R_S\cdot N:=R_S\cdot\psi(N)=\left\{\frac{x}{s}\mid x\in N, s\in
S\right\} \subseteq M_S
$$
For a subset $A\subseteq M$ we define
$A_S:=\{\frac{x}{s}\mid x\in A, s\in S\}\subseteq M_S$.
\end{notation}
\begin{rem}
Since forming the quotient module is an exact functor we can identify
the submodule $R_S\cdot N$ of $M_S$ with the quotient module $N_S$.
Also the notation $A_S$ for a set $A$ is compatible with the notation
$A_S$ if $A$ is an $R$-module.
\end{rem}
\begin{prop}
\label{up-and down}
\strut
\begin{enumerate}
\item
For each $R_S$-submodule $U$ of $M_S$ we have
$$
R_S\cdot(U\cap M)= U.
$$
\item
For each $R$-submodule $N$ of $M$ we have
$$
(R_S\cdot N)\cap M=\{x\mid x\in M, \text{ exists } s\in S
\text{ with } s\cdot x\in N \}
\supseteq N
$$
\end{enumerate}
\end{prop}
\begin{pf}
(1) Clearly $R_S\cdot(U\cap M)\subseteq U$.
Now let $\frac{x}{s}\in U$ with $s\in S$. Then
$\psi(x)=\frac{x}{1}\in U$ and
therefore $\frac{x}{t}=\frac{1}{t}\cdot\psi(x)\in U$.
\halign{&$#$\hfil\cr
(2)\quad&(N_S)\cap M&=\psi^{-1}(R_S\cdot\psi(N))\cr
&&=\{x\mid x\in M,\ \frac{x}{1}\in N_S\}\cr
&&=\{x\mid x\in M,\text{ ex. }s'\in S,\ y\in N \text{ with }
\frac{x}{1}=\frac{y}{s'} \text{ in }M_S \}\cr
&&=\{x\mid x\in M,\text{ ex. }s',t\in S,\ y\in N\text{ with }
\underbrace{ts'}_{=s\in S}\cdot x=t\cdot y\text{ in }M\}\cr
&&=\{x\mid x\in M,\text{ ex. }s\in S\text{ with }s\cdot x\in N\}
\qquad\qquad\qed\cr
}
\renewcommand{\qed}{}
\end{pf}
The operation of extending a submodule of $M$ to $M_S$ and then
restricting it back to $M$ plays an important role. Therefore an
extra name is introduced:
\begin{defn}
For any $R$-submodule $N$ of $M$ and any multiplicatively closed
subset $S$ of $R$ we define
$$
S^M(N):=(R_S\cdot N)\cap M
$$
the ``$S$-component of $N$ in $M$''.\
If no confusion can arise we will also write $S(N)$ instead of $S^M(N)$.
\end{defn}
The basic properties of these operations are summarized in the
following Proposition, the proof of which is immediate:
\begin{prop}
\label{basics-for-S-component}
Let $S,T$ be multiplicatively closed subsets of $R$; $N, \widetilde N, N_i$
submodules of $M$. Then:
\begin{enumerate}
\item
\label{RSN=RN}
$R_S\cdot S(N)=R_S\cdot N$
\item
$S(N)\supseteq N$
\item
\label{T-bigger-S}
$T\supseteq S\Longrightarrow S(T(N))=T(N)=T(S(N))$. Especially $S(S(N))=S(N)$
\item
$N\subseteq \widetilde N\Longrightarrow S(N)\subseteq S(\widetilde N)$
\item
\label{RS-ext-intersect}
$R_S\cdot(N\cap \widetilde N)=R_S\cdot N\cap R_S\cdot \widetilde N$
\item
\label{S-comp-intersect}
$S(N\cap \widetilde N)=S(N)\cap S(\widetilde N)$
\item
\label{RS-ext-union}
$\left(\bigcup\limits_i N_i\right)_S=\bigcup\limits_i (N_i)_S$
for arbitrary unions.\\
(But in general these are only sets, not modules\,!)
\item
\label{RS-ext-sum}
$R_S\cdot\sum\limits_{i} N_i=\sum\limits_{i} R_S\cdot N_i$ for
arbitrary sums.
\item
$S(N+\widetilde N)\supseteq S(N)+S(\widetilde N)$,\\
but equality does not hold in general (Example \ref{S-comp-non-add}).
\item
\label{S-comp-factor}
$S^{M/N}(0)=S^M(N)/N$
\item
If $\frak a$ is an ideal of $R$ then
$R_S\cdot(\frak a\cdot N)=(R_S\cdot\frak a)\cdot(R_S\cdot N)$
\item
\label{S-comp-M-M'}
Let $R':=R/Ann_R(M)$ and $\phi:R\rightarrow R'$
the canonical homomorphism. Then  $M$ and $N$ have  natural structures as
$R'$-modules. We denote these by $M'$ and $N'$. Let $S$ be a
multiplicatively closed subset of $R$. Then $S':=\phi(S)$ is a
multiplicatively closed subset of $R'$ and\\
${S'}^{M'}(N')=S^M(N)$ as $R$-modules and as $R'$-modules.
\end{enumerate}
\end{prop}
\begin{ex}[$S(N+\widetilde N)\nsubseteq S(N)+S(\widetilde N)$]
\label{S-comp-non-add}
\strut\\
$R:=M:=k[X,Y]$ polynomial ring in $X,Y$ over a field $k$,\\
$N:=R\cdot X$,\quad$\widetilde N:=R\cdot Y$\\
$S:=\{(X+Y)^\nu\mid \nu=0,1,2,\dots\}$.\\
Then $S(N)=N$, $S(\widetilde N)=\widetilde N$, but
$S(N+\widetilde N)=R\nsubseteq R\cdot X+R\cdot Y$.
\end{ex}
An immediate consequence of Proposition \ref{basics-for-S-component}
is the following:
\begin{prop}
The map
$$
M_S\supseteq U\mapsto U\cap M\subseteq M
$$
is an order preserving isomorphism between the lattice of the
$R_S$-submodules of $M_S$ and the lattice of those $R$-submodules of
$M$ which are $S$-components.\\
(The order is defined by ``$\,\supseteq\,$'' and the lattice operations
are ``$\,\cap\,$'',\ ``$\,+\,$'' in $M_S$ and $``\,\cap\,$'',
\ ``$\,S(\quad+\quad)\,$'' in $M$.)
\end{prop}
\begin{cor}
If $M$ is a noetherian $R$-module then $M_S$ is a noetherian $R_S$-module.
\end{cor}
There is a close connection between the primary submodules of $M$ and
$M_S$:
\begin{prop}
\label{primary-up-down}
There is a one-to-one correspondence
$$
\begin{array}{rcrclcl}
M&\supset& N&\rightmapsto&R_S\cdot N&\subset&M_S\\
M&\supset& U\cap M&\leftmapsto&U&\subset&M_S\\
\end{array}
$$
between the primary submodules $U$ of $M_S$ and those primary
submodules $N$ of $M$ whose prime ideals $\frak p$ don't intersect
with $S$, i.e. $\frak p\cap S=\emptyset$:\\
More precisely:
\begin{enumerate}
\item
If $U$ is a $\frak P$-primary submodule of $M_S$ then $U\cap M$ is a
$\frak p:=\frak P\cap R$-primary submodule of $M$, $\frak p\cap
S=\emptyset$ and $R_S\cdot (U\cap M)=U$.
\item
\label{up-P-primary}
If $N$ is a $\frak p$-primary submodule of $M$ and $\frak p\cap
S=\emptyset$ then $R_S\cdot N$ is a\\
$\frak P:=R_S\cdot\frak p$-primary
submodule of $M_S$ and $R_S\cdot N\cap M=N$.
\item
\label{S-comp-of-primary}
If $N$ is a $\frak p$-primary submodule of $M$ then:
$$
S^M(N)=\left\{
\begin{array}{ll}
N&\text{for }\frak p\cap S=\emptyset\\
M&\text{for }\frak p\cap S\ne\emptyset
\end{array}
\right.
$$
Consequently
$$
R_S\cdot N=M_S\quad\text{for }\frak p\cap S\ne\emptyset.
$$
\end{enumerate}
\end{prop}
\begin{pf}
(1) Let $U\subset M_S$ be $\frak P$-primary in $M_S$,\quad $N:=U\cap
M$ and $r\in R$ an arbitrary zero divisor for $M/N$. Then there is a
$x\in M\setminus N$ with $r\cdot x\in N$. It follows that
$\frac{r}{1}\cdot\frac{x}{1}\in R_S\cdot N=U$ (see Proposition
\ref{up-and down}). But $\frac{x}{1}\notin U$, because else
$x\in U\cap M=N$. Therefore $\frac{r}{1}$ is a zero divisor for
$M_S/U$, then $\frac{r}{1}\in\frak P$ and so
$r\in\frak P\cap R=\frak p$. Further every element $r\in\frak p$ is
nilpotent for $M$: Let $x\in M$ arbitrary. Then there is a $n\in\Bbb
N$ with $\left(\frac{r}{1}\right)^n\cdot\frac{x}{1}\in R_S\cdot N=U$
because $\frac{r}{1}\in\frak P$ is nilpotent for $M_S/U$,
and so $r^n\cdot x\in U\cap M=N$: \quad $U\cap M$ is
$\frak P\cap R$-primary. \\
(2) and (3): Let $N$ be $\frak p$-primary in $M$.\\
$1^{st}$ case: $\frak p\cap S\ne\emptyset$. Then there is an
$s\in\frak p\cap S$. $s$ is nilpotent for $M/N$, i.e. for each $x\in
M$ there is an $n\in\Bbb N$ with $s^n\cdot x\in N$ and $s^n\in S$ since $S$
is multiplicatively closed. Therefore by definition
$x\in S^M(N)$ and so $S^M(N)=M$.\\
$2^{nd}$ case: $\frak p\cap S=\emptyset$. Then $\frak P:=R_S\cdot\frak p$ is
a prime ideal of $R_S$. We claim that $R_S\cdot N$ is
$\frak P$-primary in $M_S$ and $S^M(N)=N$.\\
Proof: First we show $S(N)=N$.\quad $S(N)\supseteq N$ is always true
(Proposition \ref{basics-for-S-component}). To show the other
inclusion let $x\in S(N)$ be arbitrary. Then there is an $s\in S$ with
$s\cdot x\in N$. If $x\notin N$ then $s$ would be a zero divisor for
$M/N$ and therefore $s\in \frak p\cap S=\emptyset$~! So we have that
$S(N)=N$.\\
Now we show that $M_S/(R_S\cdot N)$ is $\frak P$-coprimary:\\
Let $\frac{r}{s}$ be an arbitrary zero divisor for $M_S/N_S$. There
is an $\frac{x}{t}\in M_S\setminus N_S$ with
$\frac{r}{s}\cdot\frac{x}{t}\in N_S$.
$\Rightarrow\frac{r\cdot x}{1}\in N_S\Rightarrow
r\cdot x\in N_S\cap M =S(N)=N$. But $x\notin N$, because else
$\frac{x}{t}\in N_S$. So $r$ is a zero divisor for $M/N$
$\Rightarrow r\in\frak p$ $\Rightarrow
\frac{r}{s}\in R_S\cdot\frak p=\frak P$.\\
Further each element of $\frak P$ is nilpotent for $M_S/N_S$ because
each element of $\frak p$ is nilpotent for $M/N$.\\
Finally: If $\frak p\cap S\ne\emptyset$ we have
$N_S=R_S\cdot N=R_S\cdot S(N)=R_S\cdot M=M_S$.
\end{pf}
\begin{prop}
\label{prim-decomp-S-comp}
Let $N$ be a proper submodule of $M$, $S$ a multiplicatively closed
subset of $R$, and
$$
N=\bigcap_{i\in I} N_i
$$
a normal (resp. primary) decomposition of $N$ in $M$, where the
$N_i$ are $\frak p_i$-primary in $M$.\\
Let
$$
I':=\{i\mid i\in I,\quad \frak p_i\cap S=\emptyset\}
$$
Then
$$
S(N)=\bigcap_{i\in I'}N_i
$$
is a normal (resp. primary) decomposition of $S(N)$ in $M$ and
$$
N_S=\bigcap_{i\in I'}(N_i)_S
$$
is a normal (resp. primary) decomposition of $N_S$ in $M_S$.
\end{prop}
\begin{pf}
\baselineskip=1.3\baselineskip
 From Proposition \ref{primary-up-down} (\ref{S-comp-of-primary}) we
know that
$$
S(N_i)=\left\{
\begin{array}{ll}
N_i&\text{ for } i\in I'\\
M               &\text{ for } i\in I\setminus I'
\end{array}
\right.
$$
It then follows from Proposition \ref{basics-for-S-component},
(\ref{S-comp-intersect}) that
$S(N)=\bigcap\limits_{i\in I'}S(N_i)\cap\bigcap\limits_{k\in I\setminus I'}
S(N_k)$\\
$=\bigcap\limits_{i\in I'}N_i$ is a primary decomposition of $S(N)$ in
$M$.\\
Clearly, if $\bigcap\limits_{i\in I}N_i$ is irredundant then
$\bigcap\limits_{i\in I'}N_i$ is irredundant too, so from a normal
decomposition of $N$ in $M$ one obtains a normal representation of
$S(N)$ in $M$.\\[1ex]
Further from Proposition \ref{basics-for-S-component}
(\ref{RS-ext-intersect}) and
Proposition \ref{primary-up-down} (\ref{up-P-primary}) and
(\ref{S-comp-of-primary})  we have\\
$ N_S=\bigcap\limits_{i\in I'}(N_i)_S\cap
\bigcap\limits_{k\in I\setminus I'}(N_k)_S=
\bigcap\limits_{i\in I'}(N_i)_S$, where the $(N_i)_S$ are
$\frak P_i:=R_S\cdot\frak p_i$-primary in $M_S$. \\
If $\frak p_i\ne \frak p_k$ then $\frak P_i\ne\frak P_k$, and if one of the
$(N_i)_S$ could be omitted in the representation
$\bigcap\limits_{i\in I'}(N_i)_S$ then the $N_i$ could be omitted in
the representation $S(N)=\bigcap\limits_{i\in I'}N_i$.\\
So, if
$N=\bigcap\limits_{i\in I}N_i$ was a normal representation then so is
$N_S=\bigcap\limits_{i\in I'} (N_i)_S$.\qed
\renewcommand{\qed}{}
\end{pf}
\subsection{Uniqueness Theorems}
\begin{prop}
Let $N=\bigcap\limits_{i\in I}N_i$ be a normal representation of $N$
in $M$,\ $N_i$ primary for $\frak p_i$ in $M$.\\
The $N_i$, whose prime ideals $\frak p_i$ are minimal in the set of
all $\frak p_i,\ i\in I$, are uniquely determined by $N$ and $M$
(i.e. they belong to any normal representation of $N$ in $M$).
\end{prop}
\begin{pf}
Let $I=\{1,\dots,r\}$ ,\ $\frak p_1$ minimal among the
$\{\frak p_1,\dots,\frak p_r\}$, and
$S:=\complement\frak p_1:=R\setminus\frak p_1$. Then for all $i\ne 1$
we have $\frak p_i\cap S\ne\emptyset$, because else
$\frak p_i\subseteq\frak p_1$ and therefore $\frak p_i=\frak p_1$
because of the minimality of $\frak p_1$, but by definition of a
normal representation $\frak p_i\ne\frak p_k$ for $i\ne k$.
Proposition \ref{prim-decomp-S-comp} yields: $S(N)=N_1$.\\
Let $N=\bigcap\limits_{j\in J}F_k$ be a second normal
decomposition of $N$ in $M$,\ $F_j$ primary for $\frak q_j$.
Again by Proposition \ref{prim-decomp-S-comp} we get
$S(N)=\bigcap\limits_{j\in J'}F_j$ with\\
$J'=\{j\mid j\in J,\ \frak q_j\cap S=\emptyset\}
=\{j\mid j\in J,\frak q_j\subseteq\frak p_1\}$.\\
$N_1$ is $\frak p_1$-primary in $M$, therefore for each $x\in M$ and
each $p\in \frak p_1$ there is a $\nu=\nu(p,x)\in\Bbb N$ with
$p^\nu\cdot x\in N_1\subseteq F_j$ for all $j\in J'$.
Then $p\in\frak q_j$, since $F_j$ is $\frak q_j$-primary in $M$, and so
finally $\frak p_1\subseteq \frak q_j$ for all  $j\in J'$.\\
On the other hand by definition of $J'$ we have $\frak q_j\subseteq\frak
p_1$, and so  $\frak q_j=\frak p_1$ for all $j\in J'$. Therefore
$J=\{j_0 \}$ contains exactly one element $j_0$,\
($\frak q_{j_0}=\frak p_1$) and therefore
$F_{j_0}=S(N)=N_1$, i.e. $N_1$ belongs also to the second normal
decomposition of $N$ in $M$.
\end{pf}
For the next Proposition we need the following facts about prime ideals:
\begin{prop}
\label{union-of-prime-ideals}
{\shape{n}\selectfont(\cite{Serre-Alg-loc}, Chap I, prop.2)}
Let $\frak a$ be an ideal and $\frak p_1,\dots\frak p_r$
finitely many prime ideals in $R$. Then
$$
\frak a\subseteq\bigcup_{i=1}^r\frak p_i
\Longleftrightarrow
\text{There exists an } i_0\in\{1,\dots r\}
\text{ with }\frak a\subseteq\frak p_{i_0}
$$
\end{prop}
\begin{cor}
\label{union-finitely-many-primes-ideal}
If $\frak p_1,\dots,\frak p_r$ are finitely many prime ideals of $R$
such that
$\frak a:=\bigcup\limits_{i=1}^r\frak p_i$ is an ideal of $R$ then
there exists an $i_0\in\{1,\dots,r\}$ with $\frak a=\frak p_{i_0}$.
\end{cor}
\begin{pf}
By  Proposition \ref{union-of-prime-ideals} there is an
$i_0\in\{1,\dots,r\}$ with
$\frak p_{i_0}\subseteq\bigcup\limits_{i=1}^r\frak p_i
=\frak a\subseteq\frak p_{i_0}$ and so $\frak a=\frak p_{i_0}$.
\end{pf}
Another immediate consequence of Proposition \ref{union-of-prime-ideals} is:
\begin{rem}
\label{a-in-a-notin-pi}
Let $\frak a$ be an arbitrary ideal, $\frak p_1,\dots,\frak p_r$ prime ideals
in $R$, and $\frak a\supsetneq \frak p_i$ for $i=1,\dots,r$.\\
Then there exists an $a\in\frak a$ with $a\notin\frak p_i$ for all
$i=1,\dots,r$.
\end{rem}
\begin{pf}
Otherwise $\frak a$ would be contained in the union of the $\frak
p_i$ and therefore in one of the $\frak p_i$.
\end{pf}
\begin{prop}
\label{complement-prime-ideal-not-in-primary-decomp}
Let $\frak p$ be a prime ideal of $R$ and
$N=\bigcap\limits_{i=1}^r N_i$ a primary decomposition of $N$ in $M$,
$N_i$ primary for $\frak p_i$.\\
Assume that $\frak p\ne\frak p_i$ for all $i=1,\dots,r$.\\
Choose an $a\in\frak p$ such that $a\notin\frak p_i$ for all
$\frak p_i$ with $\frak p_i\subseteq\frak p$\\
{\em (such an $a$ exists by Remark \ref{a-in-a-notin-pi})}\\
Let $S:=\complement\frak p$ and $T:=S\cdot\{a^\nu\mid
\nu=0,1,2,\dots\}$.\\
Then
$$
T\supsetneq S\text{\quad but\quad}  T(N)=S(N).
$$
\end{prop}
\begin{pf}
Since the $\frak p_i$ are prime ideals we have:\\
$\frak p_i\cap T=\emptyset\Longleftrightarrow
\frak p_i\cap S=\emptyset \text{ and } a\notin \frak p_i
\Longleftrightarrow \frak p_i\subseteq\frak p \text{ and }
a\notin \frak p_i$.\\
By the choice of $a$ the condition $a\notin\frak p_i$ automatically
holds for all $\frak p_i\subseteq\frak p$ and so we get:\\
$\frak p_i\cap T=\emptyset\Longleftrightarrow \frak p_i\cap S=\emptyset$.\\
Together with Proposition \ref{prim-decomp-S-comp} we obtain
$T(N)=\bigcap\limits_{\frak p_i\cap T=\emptyset}N_i =
\bigcap\limits_{\frak p_i\cap S=\emptyset}N_i =S(N)$.
On the other hand $T\supsetneq S$ since $a\in T\setminus S$.
\end{pf}
The situation is quite different for a prime ideal that occurs in a normal
decomposition:
\begin{prop}
\label{complement-prime-ideal-in-normal-decomp}
Let $N=\bigcap\limits_{i=1}^r N_i$ be a normal representation of $N$
in $M$,\\
$N_i$ primary for $\frak p_i$, and let
$\frak p\in\{\frak p_1,\dots,\frak p_r\}$,\quad $S:=\complement\frak p$.
Then:
\begin{gather*}
\text{For any multiplicatively closed subset }T\supsetneq S
\text{ of }R\\
T(N)\supsetneq S(N).
\end{gather*}
\end{prop}
\begin{pf}
Let $\frak p=\frak p_{i_0}$. Then by Proposition
\ref{prim-decomp-S-comp} $S(N)=\bigcap\limits_{i\in I'} N_i$ with
$I'=\{i\mid \frak p_i\subseteq\frak p \}\ni i_0$ is a normal
decomposition of $S(N)$.\\
Now let $T\supsetneq S$ be an arbitrary multiplicatively closed set
bigger than $S$. Then $T\cap\frak p_{i_0}\ne\emptyset$ and therefore
$T(N)=\bigcap\limits_{i\in I''}N_i$ with
$I''=\{i\mid\frak p_i\cap T=\emptyset\}\subsetneq I'$, since
$i_0\in I'\setminus I''$. It follows that $T(N)\supsetneq S(N)$
because $S(N)=\bigcap\limits_{i\in I'}N_i$ is a normal representation
and therefore $N_{i_0}$ cannot be omitted.
\end{pf}
 From Propositions \ref{complement-prime-ideal-not-in-primary-decomp}
and \ref{complement-prime-ideal-in-normal-decomp} one obtains:
\begin{cor}
\label{charact-prime-ideal-in normal-decomp}
{\shape{n}\selectfont(Compare \cite{Krull-Ringe-ohne-Endlichkeit}, Satz 12.)}
Let $\frak p$ be a prime ideal of $R$ and $N$ a proper submodule of
$M$.\\
$\frak p$ belongs to every normal representation of $N$ in $M$ if and
only if for any multiplicatively closed subset $T$ of $R$ with
$T\supsetneq\complement\frak p$ one has
$T(N)\ne\complement\frak p(N)$
(i.e. $T(N)\supsetneq\complement\frak p(N)$).
\end{cor}
Since Corollary \ref{charact-prime-ideal-in normal-decomp}
gives a characterization of the prime ideals that belong to an
arbitrary normal representation independently of that decomposition
one obtains
\begin{cor}
The set of prime ideals belonging to a normal representation of $N$ in
$M$ depends only on $N$ and $M$ and not on the representation.
\end{cor}
\section{Associated and Essential Prime Ideals}
We would like to define the ``associated'' prime ideals of a module $M$
as those that belong to a normal representation of $(0)$ in $M$. But
for non noetherian modules such a decomposition may not exist.
Nevertheless we can use the characterization given in Corollary
\ref{charact-prime-ideal-in normal-decomp} which makes sense also in
the non noetherian case (compare \cite{Krull-Ringe-ohne-Endlichkeit},
Definition on page 742):
\begin{defn}
\label{def-assprim}
Let $M$ be an arbitrary $R$-module. A prime ideal $\frak p$ of $R$ is
called an ``associated prime ideal of $M$'' iff for any
multiplicatively closed subset $T$ of $R$ with
$T\supsetneq\complement\frak p$
$$
T^M((0))\supsetneq\complement\frak p^M((0))\ .
$$
The set of all prime ideals associated to $M$ is denoted by
$\Ass(M)$.\\
If $N$ is a proper submodule of $M$ then the associated prime ideals
of $M/N$ are called the ``essential prime ideals for $N$ in $M$''.
\end{defn}
\begin{rem}
\label{ess-and-mult-closed}
$\frak p$ is essential for $N$ in $M$ iff for any multiplicatively
closed subset $T$ of $R$ with $T\supsetneq\complement\frak p$ one has
$T^M(N)\supsetneq\complement\frak p^M(N)$.
\end{rem}
\begin{pf}
By Proposition \ref{basics-for-S-component}
 (\ref{S-comp-factor}) we have\\
$T^{M/N}((0))=T^M(N)/N$ and
$\complement\frak p^{M/N}((0))=\complement\frak p^M(N)/N$.\\
Consequently $T^{M/N}((0))\supsetneq\complement\frak p^{M/N}((0))
\Longleftrightarrow T^M(N)\supsetneq\complement\frak p^M(N)$.
\end{pf}
Immediately from the definition together with Corollary
\ref{charact-prime-ideal-in normal-decomp} follows:
\begin{rem}
\label{essential-primes-for-normal-decomp}
If there exists a primary decomposition of $N$ in $M$ (e.g. if $M/N$
is noetherian) then $\frak p$ is essential for $N$ in $M$ iff $\frak
p$ belongs to a normal representation of $N$ in $M$.\\
In this case there are only {\em finitely many} essential prime
ideals for $N$ in $M$.
\end{rem}
Since in a $\frak p$-coprimary module $(0)=(0)$ is a normal
representation of $(0)$ in $M$ one has:
\begin{rem}
\label{coprimary-Ass}
If $M$ is $\frak p$-coprimary then $\Ass(M)=\{\frak p\}$\\
(The converse is also true: Corollary \ref{coprimary-iff-Assp}.)
\end{rem}
\begin{rem}\strut
\label{AssM-and-M'}
\begin{enumerate}
\item
\label{AssM-AnnM}
Each $\frak p\in\Ass(M)$ contains $\Ann_R(M)$.
\item
\label{AssM-AssM'}
Let $R':=R/\Ann_R(M)$ and $\phi:R\rightarrow R'$ the canonical
homomorphism. Then $M$ can be regarded as an $R'$-module $M'$
in a natural way. There is a one-one correspondence between $\Ass(M)$
and $\Ass(M')$, given by
$$
\Ass(M)\ni\frak p\longmapsto \phi(\frak p)\in\Ass(M')
$$
\end{enumerate}
\end{rem}
\begin{pf}
(\ref{AssM-AnnM})
If there was an $s\in\Ann_R(M)$ but $s\notin\frak p$ we would obtain
$\complement\frak p((0))=M$ and therefore also
$T((0))=M=\complement\frak p((0))$ for
all $T\supsetneq\complement\frak p$, which means that
$\frak p\notin\Ass(M)$.\\
(\ref{AssM-AssM'})
Since the prime ideals of $R'$ are in one-to-one correspondence under
$\phi$  with those prime ideals of $R$ that contains $\Ann_R(M)$, and
by (\ref{AssM-AnnM}) we know that the elements of $\Ass(M)$ contain
$\Ann_R(M)$, we only need to show that $\frak p\in \Ass(M)$ iff
$\frak p':=\phi(\frak p)\in\Ass(M')$. Obviously
$\phi\left(\complement\frak p\right)=\complement\frak p'$.
Let $T\supsetneq\complement\frak p$. Then
$T':=\phi(T)\supsetneq\complement\frak p'$ and vice versa. On the
other hand we know from \ref{basics-for-S-component} (\ref{S-comp-M-M'})
that $T^M((0))={T'}^{M'}((0))$ and also
${\complement\frak p}^M((0))={\complement\frak p'}^{M'}((0)$.
Then it follows immediately from Definition
\ref{def-assprim} that $\frak p\in\Ass(M)\Leftrightarrow\frak p'\in\Ass(M')$.
\end{pf}
\begin{prop}
\label{zero-divisors-in-Mp}
$\frak p\in\Ass(M)$ iff each element of $\frak p\cdot R_{\frak p}$ is a
zero divisor for $M_{\frak p}$.
\end{prop}
\begin{pf}
$\underline{\Rightarrow:}$ Let $\frak p\in\Ass(M)$, $p\in\frak p$,
$s\in\complement\frak p$, and
$T:=\complement\frak p\cdot\{p^\nu\mid \nu=0,1,\dots\}$.
Then $T\supsetneq\complement\frak p$ and consequently
$T((0))\supsetneq\complement\frak p((0))$, because of
$\frak p\in\Ass(M)$.
Therefore there is an $x\in M$, $x\notin\complement\frak p((0))$ with
$x\in T(0)$, i.e. there are $s'\in\complement\frak p$ and
$\nu\in\Bbb N$ such that $s'\cdot p^\nu\cdot x=0$.
($\nu>0$ since $x\notin\complement\frak p((0))$.)
Then $\frac{p^\nu}{1}\cdot\frac{x}{1}=0$ in $M_{\frak p}$ and
therefore also $(\frac{p}{s})^\nu\cdot\frac{x}{1}=0$ in
$M_{\frak p}$.  But $\frac{x}{1}\ne0$ since
$x\notin\complement\frak p(0)$ and so $(\frac{p}{s})^\nu$ and
therefore also $\frac{p}{s}$ is a zero divisor for $M_{\frak p}$.\\
$\underline{\Leftarrow:}$
Let $\frak p$ be a prime ideal of $R$ such that each element of
$\frak p\cdot R_{\frak p}$ is a zero divisor for $M_{\frak p}$,
and let $T$ be a multiplicatively closed set with
$T\supsetneq\complement\frak p$. Then there exists a $p\in T\cap \frak p$.
By assumption $\frac{p}{1}$ is a zero divisor for $M_{\frak p}$. So
there are $x\in M$ and $s\in\complement\frak p$ such that
$\frac{x}{s}\ne0$ but $\frac{p}{1}\cdot\frac{x}{s}=0$ in
$M_{\frak p}$. Then there is an $s'\in \complement\frak p$ with
$s'\cdot p\cdot x=0$ in $M$. By definition of $T$ we have $s'\cdot p\in T$
and therefore $x\in T((0))$. But $x\notin\complement\frak p$, since
otherwise $\frac{x}{s}=0$ in $M_{\frak p}$ contrary to our
assumptions. So we get $T((0))\supsetneq\complement\frak p((0))$ and
therefore $\frak p\in\Ass(M)$.
\end{pf}
Since there are no zero divisors for the zero module we get:
\begin{cor}
\label{p-Ass-Mp-not-0}
$\frak p\in\Ass(M)\Longrightarrow M_{\frak p}\ne (0)$
\end{cor}
\begin{cor}
\label{p-in Ass-is-zero-divisor}
$\frak p\in\Ass(M)\Longrightarrow$ Each element of $\frak p$ is a
zero divisor for $M$. (For the converse see Theorem
\ref{ass-zero-divisors}.)
\end{cor}
\begin{pf}
If $p\in\frak p$ then by Proposition \ref{zero-divisors-in-Mp}
$\frac{p}{1}$ is a zero divisor for $M_{\frak p}$. Therefore there is
an $x\in M$ with $\frac{x}{1}\ne0$ and $\frac{p\cdot x}{1}=0$
in $M_{\frak p}$. Then there is an $s\in\complement\frak p$ with
$p\cdot(s\cdot x)=0$ in $M$. But $s\cdot x\ne0$, since else
$\frac{x}{1}=\frac{s\cdot x}{s}=0$ in $M_{\frak p}$. So $p$ is a zero
divisor for $M$.
\end{pf}
\begin{cor}
\label{Ass-up-down}
Let $S$ be a multiplicatively closed subset of $R$ and $\frak P$ a
prime ideal of $R_S$. Then
$$
\frak P\in\Ass(M_S)\Longleftrightarrow \frak P\cap R\in\Ass(M)
$$
\end{cor}
\begin{pf}
Let $\frak p:=\frak P\cap R$. Then
$R_{\frak p}=\left(R_S\right)_{\frak P}$,\quad
$M_{\frak p}=\left(M_S\right)_{\frak P}$,\quad
$\frak p\cdot R_{\frak p}=\frak P\cdot\left(R_S\right)_{\frak P}$.\\
Therefore by Proposition \ref{zero-divisors-in-Mp}
$\frak P\in\Ass(M_S)\Longleftrightarrow$
each element of $\frak P\cdot\left(R_S\right)_{\frak P}
=\frak p\cdot R_{\frak p}$ is a zero divisor for
$\left(M_S\right)_{\frak P}=M_{\frak p}\Longleftrightarrow
\frak p\in\Ass(M)$.
\end{pf}
\begin{cor}
\label{Ass-submodule}
If $N$ is a submodule of $M$ then $\Ass(N)\subseteq\Ass(M)$.
\end{cor}
\begin{pf}
Let $\frak p\in\Ass(N)$. By Proposition \ref{zero-divisors-in-Mp}
each element of $\frak p\cdot R_{\frak p}$ is a zero divisior for
$N_{\frak p}\subseteq M_{\frak p}$ and therefore also for
$M_{\frak p}$. Again by Proposition \ref{zero-divisors-in-Mp} we get
$\frak p\in\Ass(M)$.
\end{pf}
\begin{prop}
\label{min-prime-ideal-of Ann Rpx}
Let $\frak p$ be a prime ideal of $R$, $x\in M$. The following conditions
are equivalent:
\begin{enumerate}
\item
\label{min-prime-over-ann}
$\frak p$ is minimal among the prime ideals containing $\Ann_R(x)$.
\item
\label{nilpotent-in-Rp}
$\frac{x}{1}\ne0$, and each element of $\frak p\cdot R_{\frak p}$ is
nilpotent for $R_{\frak p}\cdot x$.
\item
\label{xRp-coprimary}
$R_{\frak p}\cdot x$ is $\frak p\cdot R_{\frak p}$-coprimary.
\end{enumerate}
\end{prop}
\begin{pf}
(\ref{min-prime-over-ann}) $\Rightarrow$
(\ref{nilpotent-in-Rp}):
$\frac{x}{1}\ne0$ in $M_{\frak p}$, because else there would be a
$s\in\complement\frak p$ with $s\cdot x=0$ in $M$. Then
$s\in\Ann_R(x)\subseteq\frak p$, which contradicts
$s\in\complement\frak p$.\\
Now let $p$ be an arbitrary element of $\frak p$, and
$S:=\complement\frak p\cdot\{p^\nu\mid\nu=0,1,\dots\}$.
We show that $\Ann_{R_S}\left(\frac{x}{1}\right)=R_S$:\quad
Otherwise there would exist a prime ideal $\frak P'$ of $R_S$ with
$\frak P'\supseteq\Ann_{R_S}\left(\frac{x}{1}\right)$.
Let $\frak p':=\frak P'\cap R$. Then $\frak p'\supseteq
\Ann_{R_S}\left(\frac{x}{1}\right)\cap R\supseteq\Ann(x)$ and
$\frak p'\cap S=\emptyset$, since also
$\frak p'\cap\complement\frak p=\emptyset$, i.e.
$\frak p'\subseteq\frak p$. Now by hypothesis $\frak p$ is minimal
among the prime ideals containing $\Ann_R(x)$ and therefore
$\frak p'=\frak p$. It would follow that $p\in S\cap\frak p'$,
contradicting $\frak p'\cap S=\emptyset$. So
$\Ann_{R_S}\left(\frac{x}{1}\right)=R_S$ and consequently
$\frac{x}{1}=0$ in $M_S$.
Therefore there exist $\nu\in\Bbb N$ and $s\in\complement\frak p$
with $s\cdot p^\nu\cdot x=0$ in $M$. If follows that
$\left(\frac{p}{1}\right)^\nu\cdot\frac{x}{1}=0$ in $M_{\frak p}$.
Since $p$ was an arbitrary element of $\frak p$ we have shown that
each element of $\frac{\frak p}{1}$ and therefore also each element
of $\frak p\cdot R_{\frak p}$ is nilpotent for $R_{\frak p}\cdot x$.
\\[1ex]
(\ref{nilpotent-in-Rp}) $\Rightarrow$
(\ref{xRp-coprimary}):
By assumption each element of $\frak p\cdot R_{\frak p}$ is nilpotent
for $R_{\frak p}\cdot x$ and therefore a zero divisor for
$R_{\frak p}\cdot x$, because $R_{\frak p}\cdot x\ne(0)$.
Since the elements of $R_{\frak p}\setminus\frak p\cdot R_{\frak p}$
are the units of $R_{\frak p}$, $\frak p\cdot R_{\frak p}$ is the set
of all zero divisors for $R_{\frak p}\cdot x$. Consequently
$R_{\frak p}\cdot x$ is $\frak p\cdot R_{\frak p}$-coprimary.
\\[1ex]
(\ref{xRp-coprimary}) $\Rightarrow$
(\ref{min-prime-over-ann}):
$\frak p\supseteq\Ann(x)$, because else there would exist an
$s\in R\setminus\frak p$ with $s\cdot x=0$ and therefore
$\frac{x}{1}=0$ in $M_{\frak p}$ which cannot be since
$R_{\frak p}\cdot x$ is coprimary and therefore $\ne (0)$.\\
$\frak p$ is minimal among the prime ideals containing $\Ann(x)$:
Let $\frak p'$ be a prime ideal with
$\frak p\supseteq\frak p'\supseteq\Ann(x)$ and $p\in\frak p$. Then
by assumption $\frac{p}{1}$ is nilpotent for $R_{\frak p}\cdot x$.
Then there exist $s\in R\setminus\frak p$ and $\nu\in\Bbb N$ with
$s\cdot p^\nu\in\Ann(x)\subseteq\frak p'$. But
$s\notin\frak p\supseteq\frak p'$. It follows that
$p\in\frak p'$ and therefore $\frak p=\frak p'$.
\end{pf}
\begin{cor}
\label{min-primes-Ann-Ass}
Let $0\ne x\in M$ be an arbitrary element of $M$ and $\frak p$ minimal
among the prime ideals containing $\Ann(x)$. Then $\frak p\in\Ass(M)$.
\end{cor}
\begin{pf}
By Proposition \ref{min-prime-ideal-of Ann Rpx} $\frac{x}{1}\ne0$
in $M_{\frak p}$ and each element of $\frak p\cdot R_{\frak p}$
is nilpotent for $R_{\frak p}\cdot x$ and therefore a zero divisor
for $R_{\frak p}\cdot x$ and hence also for $M_{\frak p}$. From
Proposition \ref{zero-divisors-in-Mp}
then follows that $\frak p\in\Ass(M)$.
\end{pf}
\begin{defn}
\label{def-Ass-essential-first-kind}
\strut
\begin{enumerate}
\item
\label{def-Ass1}
A prime ideal $\frak p$ of $R$ is called
``associated of the first kind to $M$''
iff there exists an $x\in M$ such that $\frak p$ is minimal among all
prime ideals that contain $\Ann_R(x)$:
$$
\Ass_1(M):=\{\frak p\mid \frak p \text{ associated of the first
kind to $M$}\}
$$
\item
\label{def-essential1}
A prime ideal $\frak p$ of $R$ is called ``essential of the first
kind for $N$ in $M$'' iff $\frak p$ is associated of the first kind to $M/N$ .
\end{enumerate}
\end{defn}
\begin{rem}
\label{Ass1-rems}
\strut
\begin{enumerate}
\item
In \cite{Bourbaki-Alg-Comm-4} \S 1 exercise 17 the prime ideals
which we call associated of the first kind to $M$ are called ``faiblement
associ\'e \`a $M$''.
\item
\label{Ass1-in-Ass}
$\Ass_1(M)\subseteq\Ass(M)$, but in general equality does not hold
(Example \ref{Ass-not-Ass1}).
\item
\label{M-not-zero-Ass1}
$M\ne(0)\Longleftrightarrow\Ass_1(M)\ne\emptyset
\Longleftrightarrow\Ass(M)\ne\emptyset$.
\item
\label{essential-primes-for-ideal}
Let $\frak a$ be a proper ideal in $R$.
The prime ideals which are minimal in the set of all prime ideals
containing $\frak a$ are essential of the first kind for $\frak a$ in $R$.
\item
\label{Ass1M-Ass1M'}
Let $R':=R/\Ann_R(M)$ and $\phi:R\rightarrow R'$ the canonical
homomorphism. Then $M$ can be regarded as an $R'$-module $M'$
in a natural way. There is a one-to-one correspondence between $\Ass_1(M)$
and $\Ass_1(M')$, given by
$$
\Ass_1(M)\ni\frak p\longmapsto \phi(\frak p)\in\Ass_1(M')
$$
\end{enumerate}
\end{rem}
\begin{pf}
(\ref{Ass1-in-Ass}) follows from Corollary \ref{min-primes-Ann-Ass}.
\\[1ex]
(\ref{M-not-zero-Ass1})
$M\ne(0)\Rightarrow$ Ex. $0\ne x\in M\Rightarrow\Ann_R(x)\ne R
\Rightarrow$ Ex. prime ideal $\frak p\supseteq\Ann_R(x)$ and
therefore there also exists a prime ideal $\frak p'$ which is minimal
among the prime ideals containing $\Ann_R(x)$. By definition $\frak
p'\in\Ass_1(M)\subseteq\Ass(M)$.
\\
Conversely: If $\Ass(M)\ne\emptyset$ then $M\ne(0)$ (and then also
$\Ass_1(M)\ne\emptyset$ as was just shown), since by
Corollary \ref{p-Ass-Mp-not-0}
one even has $M_{\frak p}\ne(0)$ for all $\frak p\in\Ass(M)$.
\\[1ex]
(\ref{essential-primes-for-ideal})
For the residue class $\bar1\in R/\frak a$ we have
$\Ann_R(\bar1)=\frak a$. Therefore the prime ideals $\frak p$ which
are minimal among the prime ideals containing $\frak a$  belong to
$\Ass_1(R/\frak a)$, which means that $\frak p$ is essential of the
first kind for $\frak a$ in $R$.
\\[1ex]
(\ref{Ass1M-Ass1M'})
This follows immediately from the fact that for each $x\in M$ we have\\
$\Ann_{R'}(\phi(x))=\phi\left(\Ann_R(x)\right)$. (See also the proof of
Remark \ref{AssM-and-M'} (\ref{AssM-AssM'}) .)
\end{pf}
With respect to quotient modules the elements of $\Ass_1(M)$ behave
similar to those of $\Ass(M)$ (Corollary \ref{Ass-up-down}):
\begin{prop}
\label{Ass1-up-down}
Let $S$ be a multiplicatively closed subset of $R$ and $\frak P$ a
prime ideal of $R_S$. Then
$$
\frak P\in\Ass_1(M_S)\Longleftrightarrow \frak P\cap R\in\Ass_1(M)
$$
\end{prop}
\begin{pf}
We use the same notations as in the proof of Corollary
\ref{Ass-up-down}.\\
$\frak P\in\Ass_1(M)\Leftrightarrow$ there
exists $x\in M$ such that $\frak P$ is minimal among the prime ideals
containing $\Ann_{R_S}\left(\frac{x}{1}\right)$. By propositon
\ref{min-prime-ideal-of Ann Rpx} this means that $\frac{x}{1}\ne0$
and each element of $\frak P\cdot\left(R_S\right)_{\frak P}=
\frak p\cdot R_{\frak p}$ is nilpotent for
$\left(R_S\right)_{\frak P}\cdot\frac{x}{1}=R_{\frak p}\cdot x$,
which again by Proposition \ref{min-prime-ideal-of Ann Rpx},
is equivalent to $\frak p$ being minimal among the prime ideals of
$R$ containing $\Ann_R(x)$, i.e. $\frak p\in\Ass_1(M)$.
\end{pf}
\begin{prop}
\label{union-of-Ass1-primes}
\begin{enumerate}
\item
\label{each-p-is-union}
Each $\frak p\in\Ass(M)$ is the union of certain\\
$\frak p'\in\Ass_1(M)$.
More exactly:
$$
\frak p=\bigcup_{\frak P'\in\Ass_1(M_{\frak p})}\frak P'\cap R
$$
\item
\label{finite-union-of Ass1}
If $\frak p'_1,\dots,\frak p'_r\in\Ass_1(M)$ are {\em finitely many} prime
ideals of $\Ass_1(M)$ such that
$\bigcup\limits_i\frak p'_i=:\frak p$ is a prime
ideal, then $\frak p$ is equal to one of the $\frak p'_i$ and
therefore $\frak p\in\Ass_1(M)$.
\item
\label{if-union-is-prime}
If $\frak p_i\in\Ass(M)$ and $\bigcup\limits_i\frak p_i=:\frak p$ is a prime
ideal then $\frak p\in\Ass(M)$.
\end{enumerate}
\end{prop}
\begin{pf}
(\ref{each-p-is-union})
Let $\frak p\in\Ass(M)$ and
$\frak P'\in\Ass_1\left(M_{\frak p}\right)$.
Then $\frak P'\cap R\subseteq\frak p$,
and by Proposition \ref{Ass1-up-down}
$\frak P'\cap R\in\Ass_1(M)$. Therefore
$\frak p\supseteq
\bigcup\limits_{\frak P'\in\Ass_1\left(M_{\frak p}\right)}
\frak P'\cap R$.\\
Conversely: Let $p\in\frak p$. Then by Proposition
\ref{zero-divisors-in-Mp}\quad $\frac{p}{1}$ is a zero divisor for
$M_{\frak p}$. Therefore there exists
$0\ne\frac{x}{1}\in M_{\frak p}$ with
$\frac{p}{1}\in
\Ann_{R_{\frak p}}\left(R_{\frak p}\cdot\frac{x}{1}\right)$.
Let $\frak P'_x$ be minimal among the prime ideals containing
$\Ann_{R_{\frak p}}\left(R_{\frak p}\cdot\frac{x}{1}\right)$. Then by
Definition \ref{def-Ass-essential-first-kind}\ \
$\frak P'_x\in\Ass_1(M_{\frak p})$ and
$p\in\frak P'_x\cap R$. Therefore $\frak p\subseteq
\bigcup\limits_{\frak P'\in\Ass_1\left(M_{\frak p}\right)}
\frak P'\cap R$.
\\[1ex]
(\ref{finite-union-of Ass1})
This is an immediate consequence of
Corollary \ref{union-finitely-many-primes-ideal}.
\\[1ex]
(\ref{if-union-is-prime})
Let $\frak p=\bigcup\limits_i\frak p_i$ with $\frak p_i\in\Ass(M)$ be
a prime ideal. Then $\frak p\supseteq\frak p_i$ and therefore
$\frak P_i:=R_{\frak p}\cdot\frak p_i$ is a prime ideal of
$R_{\frak p}$ and
$M_{\frak p_i}=\left(M_{\frak p}\right)_{\frak P_i}$.
Now let $p\in\frak p$ be an arbitrary element. We show that
$\frac{p}{1}$ is a zero divisor of $M_{\frak p}$:\\
Since $p\in\frak p_{i_0}$ for some $i_0$ and
$\frak p_{i_0}\in\Ass(M)$ it follows from Proposition
\ref{zero-divisors-in-Mp} that
$\frac{p}{1}\in R_{\frak p_{i_0}}$ is a zero divisor for
$M_{\frak p_{i_0}}$. Then there is an $x\in M$ with
$\frac{x}{1}\ne0$ in $M_{\frak p_{i_0}}$ but
$\frac{p}{1}\cdot\frac{x}{1}=0$ in $M_{\frak p_{i_0}}$.
Then there is a $\rho\in R_{\frak p}\setminus\frak P_{i_0}$
with $\rho\cdot\frac{p}{1}\cdot\frac{x}{1}=0$ in $M_{\frak p}$,
since $M_{\frak p_{i_0}}=\left(M_{\frak p}\right)_{\frak P_{i_0}}$.
But $\rho\cdot\frac{x}{1}\ne0$ in $M_{\frak p}$, because else
$\frac{x}{1}=0$ in
$M_{\frak p_{i_0}}=\left(M_{\frak p}\right)_{\frak P_{i_0}}$.
Consequently $\frac{p}{1}\in R_{\frak p}$ is a zero divisor for
$M_{\frak p}$. Proposition \ref{zero-divisors-in-Mp}
then yields $\frak p\in\Ass(M)$.
\end{pf}
An immediate consequence of Proposition \ref{union-of-Ass1-primes}
(\ref{each-p-is-union}) and (\ref{finite-union-of Ass1}) is:
\newpage
\begin{cor}\strut
\begin{enumerate}
\item
Each $\frak p\in\Ass(M)$ contains a\/ $\frak p'\in\Ass_1(M)$. Therefore
all minimal elements of $\Ass(M)$ (if there are any) lie in $\Ass_1(M)$.
\item
If $\Ass_1(M)$ is a finite set then $\Ass(M)=\Ass_1(M)$.
(E.g. if $M$ is noetherian or if $(0)$ has a primary decomposition.)
\end{enumerate}
\end{cor}
\begin{thm}
\label{ass-zero-divisors}
$$
\bigcup\limits_{\frak p'\in\Ass_1(M)}\frak p'=
\bigcup\limits_{\frak p\in\Ass(M)}\frak p=
\{r\mid r\in R,\ r\text{ zero divisor for }M \}.
$$
\end{thm}
\begin{pf}
 From $\Ass_1(M)\subseteq\Ass(M)$ follows
$\bigcup\limits_{\frak p'\in\Ass_1(M)}\frak p'\subseteq
\bigcup\limits_{\frak p\in\Ass(M)}\frak p$, and because of
Proposition \ref{union-of-Ass1-primes} (\ref{each-p-is-union})
one has also the converse inclusion. That shows the first equality.\\
To show the second equality we only need to show that every zero divisor for
$M$ lies in a $\frak p\in\Ass(M)$, because the other inclusion
follows from Corollary \ref{p-in Ass-is-zero-divisor}.
Let $r\in R$ be a zero divisor for $M$. Then there exists an
$x\in M$, $x\ne0$ with $r\in\Ann_R(x)$. Since $\Ann_R(x)\ne R$
there exist prime ideals in $R$ containing $\Ann_R(x)$. A minimal
element among these primes belongs to $\Ass(M)$ by Corollary
\ref{min-primes-Ann-Ass} and contains $r$.
\end{pf}
Since $\Ass(M)$ describes the zero divisors for $M$
it is plausible that there is also connection with the annihilators
of submodules of $N$:
\begin{prop}
\label{ess-Ann-Ass-M-prim-decomp-0}
Let $N\ne(0)$ be a finitely generated submodule of $M$, then
\begin{enumerate}
\item
\label{ess-Ann-N-Ass-M}
The essential prime ideals for $\Ann_R(N)$ in $R$ belong to $\Ass(M)$
\item
\label{ess1-Ann-N-Ass1-M}
The essential prime ideals of the first kind for $\Ann_R(N)$ in $R$\\
belong to $\Ass_1(M)$
\item
\label{primary-decomp-0-decomp-Ann}
If there exists a primary decomposition of $(0)$ in $M$ then there is
also a primary decomposition of $\Ann_R(N)$ in $R$:
\end{enumerate}
\end{prop}
\begin{pf}
(\ref{ess-Ann-N-Ass-M}):
Let $\frak p$ be essential for $\Ann_R(N)$ in $R$ and $T$ a
multiplicatively closed subset of $R$ with
$T\supsetneq\complement\frak p$.
Then by Remark \ref{ess-and-mult-closed}\ \
$T^R\left(\Ann_R(N)\right)\supsetneq
{\complement\frak p}^R\left(\Ann_R(N)\right)$.
We will show that $R_{\frak p}\cdot T^R\left(\Ann_R(N)\right)
\supsetneq
R_{\frak p}\cdot{\complement\frak p}^R\left(\Ann_R(N)\right)$:\\
Trivially ``$\supseteq$'' holds. To show the inequality first remark
that by Proposition \ref{basics-for-S-component} (\ref{T-bigger-S})
$T^R\left(\Ann_R(N)\right)=
{\complement\frak p}^R\left(T^R(\Ann_R(N))\right)$ is a
$\complement\frak p$-component
because of $T\supseteq\complement\frak p$.
Therefore from
$R_{\frak p}\cdot T^R\left(\Ann_R(N)\right)
=
R_{\frak p}\cdot{\complement\frak p}^R\left(\Ann_R(N)\right)$
by taking the inverse images in $R$ one would obtain
(using Proposition \ref{basics-for-S-component} (\ref{RSN=RN}) and the
definition of the $\complement\frak p$-components):
$T^R\left(\Ann_R(N)\right)=
R_{\frak p}\cdot T^R\left(\Ann_R(N)\right)\cap R=
R_{\frak p}\cdot{\complement\frak p}^R\left(\Ann_R(N)\right)\cap R=
{\complement\frak p}^R\left(\Ann_R(N)\right)$, contradiction.\\
Further, again by Proposition \ref{basics-for-S-component}
(\ref{RSN=RN}), we know that
$R_{\frak p}\cdot\complement\frak p\left(\Ann_R(N)\right)=
R_{\frak p}\cdot\Ann_R(N)$, and
$R_{\frak p}\cdot\Ann_R(N)=\Ann_{R_{\frak p}}(N_{\frak p})$ , because
$N$ is finitely generated
(Remark \ref{transporteurs} (\ref{transp-quotmodule})), so that finally
$R_{\frak p}\cdot T^R\left(\Ann_R(N)\right)\supsetneq
\Ann_{R_{\frak p}}(N_{\frak p})$.\\
Let $r\in T^R\left(\Ann_R(N)\right)$ with $\frac{r}{1}\in
R_{\frak p}\cdot T^R\left(\Ann_R(N)\right)\setminus
\Ann_{R_{\frak p}}(N_{\frak p})$, i.e. there is an $x\in N$ with
$\frac{r}{1}\cdot\frac{x}{1}\ne(0)$ in
$N_{\frak p}\subseteq M_{\frak p}$ and so
$r\cdot x\notin{\complement\frak p}^M(0)$.
But $\frac{r}{1}\in R_T\cdot T^R\left(\Ann_R(N)\right)=
R_T\cdot\Ann_R(N)=\Ann_{R_T}(N_T)$ since $N$ is finitely generated
(Remark \ref{transporteurs} (\ref{transp-quotmodule})).
Therefore $\frac{r}{1}\cdot\frac{x}{1}=0$ in $N_T\subseteq M_T$ and so
$r\cdot x\in T^M((0))$.\\
So we have shown that for any
$T\supsetneq\complement\frak p$ we have
$T^M((0))\supsetneq{\complement\frak p}^M((0))$ and consequently
$\frak p\in\Ass(M)$.
\\[1ex]
(\ref{ess1-Ann-N-Ass1-M}): Now let $\frak p$ be essential of the
first kind for $\Ann_R(N)$ in $R$,\\
i.e. $\frak p\in\Ass_1(R/\Ann_R(N))$. By definition there is an $r\in R$
such that $\frak p$ is minimal among the prime ideals containing the
annihilator $\Ann_R(\bar r)$ of the residue class $\bar r$ of $r$ mod
$\Ann_R(N)$. Then by Proposition \ref{min-prime-ideal-of Ann Rpx}
\quad $\frac{\bar r}{1}\ne 0$, and for each
$\frac{p}{1}\in\frak p\cdot R_{\frak p}$ there exists a
$\nu\in\Bbb N$ with $\frac{p^\nu\cdot r}{s^\nu}\in
\left(\Ann_R(N)\right)_{\frak p}=\Ann_{R_{\frak p}}(N_{\frak p})$,
since $N$ is finitely generated
(Remark \ref{transporteurs} (\ref{transp-quotmodule})).
Then there is an $x\in N$ with
$\frac{r\cdot x}{1}\ne0$ in $N_{\frak p}\subseteq M_{\frak p}$, but
for each $\frac{p}{s}\in\frak p\cdot R_{\frak p}$ there is a
$\nu\in\Bbb N$ with
$\left(\frac{p}{s}\right)^\nu\cdot\frac{r\cdot x}{1}=0$ in
$M_{\frak p}$, showing that each element of $\frak p\cdot R_{\frak p}$
is nilpotent for $R_{\frak p}\cdot(r\cdot x)$. Proposition
\ref{min-prime-ideal-of Ann Rpx} then yields that $\frak p$ is minimal
among the prime ideals containing $\Ann_R(r\cdot x)$ and therefore
$\frak p\in\Ass_1(M)$.
\\[1ex]
(\ref{primary-decomp-0-decomp-Ann}):
Assume that there is a primary decomposition
$(0)=\bigcap\limits_{i=1}^n F_i$, with $F_i$ primary in $M$. Then
by Remark \ref{transporteurs}
$\Ann_R(N)=\left((0):N\right)=\bigcap\limits_{i=1}^n\left(F_i:N\right)$
and\\
$
(F_i:N)=
\left\{\begin{array}{ll}
\text{ primary in } R&\text{ if }F_i\nsupseteq N\\
R                    &\text{ if }F_i\supseteq N
\end{array}\right.,
$
because $N$ is finitely generated.
Hence $\Ann_R(N)$ has a primary decomposition in $R$.
\end{pf}
\begin{rem} If $N$ is not finitely generated
Proposition \ref{ess-Ann-Ass-M-prim-decomp-0} is not true
as is shown in Example \ref{0-primary-decomp-Ann-not ass}
\end{rem}
In the proof of Proposition \ref{ess-Ann-Ass-M-prim-decomp-0} we used
the following
\begin{defn}
\label{transpdefs}
Let $N$ and $U$ be subsets of the $R$-module $M$
$$
(N:U):=\{r\mid r\in R,\ r\cdot U\subseteq N\}
$$
\end{defn}
\begin{rem}
\label{transporteurs}
Let $N$ and $U$ be subsets of $M$.
\begin{enumerate}
\item
\label{submodule-ideal}
If $N$ is a $R$-submodule of $M$ then $(N:U)$ is an ideal of $R$
and\\
$(N:U)=(N:\langle U\rangle)$, where $\langle U\rangle$ denotes the
$R$-module generated by $U$.
\item
\label{Ann-transp}
$\Ann_R(U)=\left((0):U\right)$
\item
\label{intersect-transp}
For arbitrary intersections we have
$$
\left(\left(\bigcap\limits_{i\in I}F_i\right):U\right)=
\bigcap\limits_{i\in I}\left(F_i:U\right)
$$
\item
\label{primary-transp}
Let $F$ be $\frak p$-primary in $M$ and $U$ a finitely generated
submodule of $M$. Then
$$
(F:U)\text{ is }\left\{
\begin{array}{ll}
\frak p\text{-primary in }R &\text{ if }F\nsupseteq U\\
=R &\text{ if }F\supseteq U
\end{array}
\right.
$$
\item
\label{transp-quotmodule}
Let $S$ be a multiplicatively closed subset of $R$,\\
$N$, $U$ submodules of $M$. Then
$$
\left(N_S:U_S\right)\supseteq R_S\cdot\left(N:U\right)\ .
$$
If $U$ is finitely generated equality holds.\\
Especially: If $U$ is finitely generated then
$\Ann_{R_S}(U_S)=R_S\cdot\Ann_R(U)$
\end{enumerate}
\end{rem}
\begin{pf}
(\ref{submodule-ideal}),
(\ref{Ann-transp}),
(\ref{intersect-transp})
are obvious.
\\[1ex]
(\ref{primary-transp}):
Trivially $(F:U)=R$ if $F\supseteq U$.\\
Let $F\nsupseteq U$. Then $1\notin(F:U)$, i.e. $(F:U)\subsetneq R$.
By hypothesis $F$ is $\frak p$-primary in $M$. Then for each
$p\in\frak p$ and each $x\in U$ there is a $\nu\in\Bbb N$ with
$p^\nu\cdot x\in F$. Now by hypothesis $U$ is finitely generated. Given a
$p\in\frak p$ we define $n$ to the maximum of the $\nu_i$ such that for a
finite set of generators $x_i$ of $U$ we have
$p^{\nu_i}\cdot x_i\in F$. Then  by (\ref{submodule-ideal})\
$p^n\cdot U\subseteq F$, which means $p^n\cdot R\subseteq(F:U)$. It
follows that each element of $\frak p$ is nilpotent for
$R/(F:U)$. On the other hand, if $r\in R$ is any zero divisor for
$R/(F:U)$ there exists a $y\in R\setminus(F:U)$ with $r\cdot y\in (F:U)$,
i.e. $y\cdot U\nsubseteq F$, but $r\cdot y\cdot U\subseteq F$.
Therefore there is a $z\in y\cdot U$, $z\notin F$ with
$r\cdot z\in r\cdot y\cdot U\subseteq F$. So $r$ is also a zero divisor for
$M/F$ and hence $r\in\frak p$ since $F$ is $\frak p$-primary in $M$,
and we have shown that $(F:U)$ is $\frak p$-primary in $R$.
\\[1ex]
(\ref{transp-quotmodule}):
 From $r\cdot U\subseteq N$ follows
$\frac{r}{s}\cdot U_S\subseteq N_S$ for all $s\in S$, and therefore\linebreak
$R_S\cdot(N:U)\subseteq\left(N_S:U_S\right)$.
\\[1ex]
Now let $U$ be finitely generated, $U=\langle x_1,\dots,x_r\rangle$, and
let $\frac{r}{s}\in\left(N_S:U_S\right)$ be arbitrary.
Then for $i=1,\dots r$\quad $\frac{r}{s}\cdot\frac{x_i}{1}\in N_S$,
i.e. for all $i=1,\dots r$ there exists an $s_i\in S$ with
$s_i\cdot r\cdot x_i\in N$. Define $s':=s_1\cdots s_r$. Then
$s'\cdot r\cdot x_i\in N$ for all $i$ and therefore
$s'\cdot r\in (N:U)$. It follows
$\frac{r}{1}\in R_S\cdot(N:U)$ and therefore also
$\frac{r}{s}\in R_S\cdot(N:U)$.
\end{pf}
\begin{rem}
Without the assumption that $U$ is finitely generated the conclusion
of Remark \ref{transporteurs} (\ref{primary-transp}) may be false
as can be seen from Example \ref{0-primary-decomp-Ann-not ass} with
$F:=(0)$, $U:=M$. Then $(0)$ is $\frak p$-primary in $M$, but
$(F:U)=\Ann_R(M)=(0)$ is not $\frak p$-primary in $R$.
\end{rem}
\vspace{2ex}\noindent
In the classical case of noetherian modules
one defines the associated prime ideals of $M$
as those prime ideals which are annihilators of elements of $M$ and
not just minimal elements in the set of all prime ideals containing
the annihilator of an element. Here we will denote the set of these
prime ideals by $\Ass_0(M)$:
\begin{defn}
$$
\Ass_0(M):=\{\frak p\mid \frak p
{\size{10}{12pt}\selectfont\text{ prime ideal of $R$ such that
there is an $x\in R$ with }}\frak p=\Ann_R(x)\}.
$$
\end{defn}
\begin{rem}
\label{Ass0-rems}
\strut
\begin{enumerate}
\item
Clearly by definition $\Ass_0(M)\subseteq\Ass_1(M)$. But in general
equality does not hold. (See Example \ref{Ass=Ass1-not-Ass0}.)
\item
If $M$ is noetherian, then $\Ass(M)=\Ass_1(M)=\Ass_0(M)$.
(See Theorem \ref{Ass-M-noetherian}.)
In the non noetherian case $\Ass_0(M)$ is
not very useful. For instance it may happen that
$\Ass_0(M)=\emptyset$ although there exists a primary decomposition
of $(0)$ in $M$. (See Example \ref{Ass=Ass1-not-Ass0}.)
\item
\label{Ass0M-Ass0M'}
Let $R':=R/\Ann_R(M)$ and $\phi:R\rightarrow R'$ the canonical
homomorphism. Then $M$ can be regarded as an $R'$-module $M'$
in a natural way. There is a one-to-one correspondence between $\Ass_0(M)$
and $\Ass_0(M')$, given by\hfil
$\Ass_0(M)\ni\frak p\longmapsto \phi(\frak p)\in\Ass_0(M')$\\
(See the proof of Remark \ref{Ass1-rems} (\ref{Ass1M-Ass1M'}) .)
\end{enumerate}
\end{rem}
\begin{thm}
\label{Ass-M-noetherian}
Let $M$ be an $R$-module and $\frak p$ a prime ideal of $R$.
\begin{enumerate}
\item
\label{p-coprim-p-Ass}
If there exists a $\frak p$-coprimary submodule $U$ of $M$, then
$\frak p\in\Ass_1(M)$.
\item
\label{0-primary-decomp-cyclic-p-primary}
If $\frak p\in\Ass(M)$ and $(0)$ has a primary decomposition in $M$
then there exists a cyclic $\frak p$-coprimary submodule $U$ of $M$.
\item
\label{0-primary-decomp-cyclic-p-primary-finite-p}
If there exists a primary decomposition of $(0)$ in $M$ and $\frak p$
is finitely generated then:
$\frak p\in \Ass(M)\Longleftrightarrow \frak p\in\Ass_0(M)$.
Especially, if $M$ and $R$ are both noetherian then
$\Ass(M)=\Ass_1(M)=\Ass_0(M)$
\item
\label{R-noetherian-Ass1=Ass0}
If $R$ is noetherian then $\Ass_1(M)=\Ass_0(M)$
\item
\label{M-noetherian-Ass=Ass0}
If $M$ is noetherian then $\Ass(M)=\Ass_1(M)=\Ass_0(M)$
\end{enumerate}
\end{thm}
\begin{pf}
(\ref{p-coprim-p-Ass}) Let $U$ be a $\frak p$-coprimary submodule of
$M$. Then by definition $U\ne(0)$. Choose an $0\ne x\in U$. By
Remark \ref{submodule-coprimary-is coprimary} $R\cdot x$ is
$\frak p$-coprimary and therefore by
Proposition \ref{primary-up-down}
$R_{\frak p}\cdot x$ is $\frak p\cdot R_{\frak p}$-coprimary. By
Proposition \ref{min-prime-ideal-of Ann Rpx}
and the definition of $\Ass_1(M)$ we get $\frak p\in
\Ass_1(M)$.
\\[1ex]
(\ref{0-primary-decomp-cyclic-p-primary}) By hypothesis there exists
a primary and therefore also a normal decomposition of $(0)$ in
$M$:\\
$(0)=\bigcap\limits_{i=1}^n F_i$ with $F_i$ being $\frak p_i$-primary
in $M$. By Remark \ref{essential-primes-for-normal-decomp}
$\Ass(M)=\{\frak p_1,\dots,\frak p_n\}$.
Without loss of generality we may assume that $\frak p=\frak p_1$.
We distinguish two cases:\\
1$^{st}$ case: $n=1$. Then $(0)=F_1$ is $\frak p$-primary in $M$, i.e. $M$
is $\frak p$-coprimary, especially $M\ne(0)$. Let $0\ne x\in M$. Then
by Remark \ref{submodule-coprimary-is coprimary}
$U:=R\cdot x$ is a cyclic $\frak p$-coprimary submodule of $M$.\\
2$^{nd}$ case: $n\ge 2$. Let
$\widetilde U:=\bigcap\limits_{i=2}^n F_i$ Then $\widetilde U\ne(0)$, since
$(0)=F_1\cap\dots\cap F_n$ is reduced. We show that $\widetilde U$ is
$\frak p$-coprimary:\\
Since $F_1$ is $\frak p$-primary in $M$, for each $p\in\frak p$ and
each $x\in M$ there exists a $\nu\in\Bbb N$ with
$p^\nu\cdot x\in F_1$. Applying this to an $x\in\widetilde U$ we get
$p^\nu\cdot x\in F_1\cap\widetilde U=(0)$. Therefore each element of
$\frak p$ is nilpotent for $\widetilde U$.\\
Now let $r$ be an arbitrary zero divisor for $\widetilde U$.
Then there is an $0\ne x\in\widetilde U$ with $r\cdot x=0$. We have
$x\notin F_1$, since else $x\in F_1\cap\widetilde U=(0)$, but
$r\cdot x=0\in F_1$. Therefore $r$ is a zero divisor for $M/F_1$,
hence $r\in\frak p$. It follows that $\widetilde U$ is
$\frak p$-coprimary. Then any $0\ne x\in\widetilde U$ generates a
cyclic $\frak p$-coprimary submodule of $M$.
\\[1ex]
(\ref{0-primary-decomp-cyclic-p-primary-finite-p})
Since $\Ass(M)\supseteq\Ass_1(M)\supseteq\Ass_0(M)$
all we have to show is that for each
$\frak p\in\Ass(M)$ there exists a cyclic submodule of $M$ whose
annihilator is $\frak p$.
By (\ref{0-primary-decomp-cyclic-p-primary}) we have a cyclic
$\frak p$-coprimary submodule $U$ of $M$.
Since each element of $\frak p$ is nilpotent for $U$ and
$\frak p$ and $U$ are finitely generated there exists a
$\nu\in\Bbb N$ with $\frak p^\nu\cdot U=(0)$. Choose $\nu$ minimal
with that property, then $\frak p^{\nu-1}\cdot U\ne(0)$. Let
$0\ne y\in\frak p^{\nu-1}\cdot U$. Then $\frak p\cdot y=0$ and
every $r\in R$ with $r\cdot y=0$ lies in $\frak p$ since
$R\cdot y$ is $\frak p$-coprimary as a submodule of $U$. Therefore
$\Ann(y)=\frak p$.
\\[1ex]
(\ref{R-noetherian-Ass1=Ass0})
Let $R$ be noetherian. Since $\Ass_1(M)\supseteq\Ass_0(M)$ all we have
to show is that each $\frak p\in\Ass_1(M)$ belongs to $\Ass_0(M)$:\\
For $\frak p\in\Ass_1(M)$ there exists an $x\in M$ such that $\frak p$
is minimal among the prime ideals containing $\frak a:=\Ann_R(x)$.
Let $V:=R\cdot x\subseteq M$. Then by definition $\frak p\in\Ass_1(V)$.
But $V\cong R/\frak a$ is a noetherian $R$-module since $R$ is
noetherian. By~(\ref{0-primary-decomp-cyclic-p-primary-finite-p})
we then get $\frak p\in\Ass_0(V)$, i.e. there is an element
$y\in V$ with $\frak p=\Ann_R(y)$. But since $y\in V\subseteq M$ we
obtain $\frak p\in\Ass_0(M)$.
\\[1ex]
(\ref{M-noetherian-Ass=Ass0})
Let $M$ be noetherian, $\frak a:=\Ann_R(M)$, and $R':=R/\frak a$.
Then $M$ can be regarded in a natural way as an $R'$-module $M'$ and by
\cite{Nagata-Local-Rings} Corollary (3.17) $R'$ is a noetherian ring.
Therefore by (\ref{0-primary-decomp-cyclic-p-primary-finite-p}) we
have $\Ass(M')=\Ass_1(M')=\Ass_0(M')$ as an $R'$-module.
But then the same equality holds for $M$ as an $R$-module, because of
the one-to-one correspondence between the respective $\Ass$. (See
Remarks \ref{AssM-and-M'}, \ref{Ass1-rems}, \ref{Ass0-rems}.)
\end{pf}
\begin{cor}
\label{Ass12-exactq-sequ}
Let
$0\rightarrow N\rightarrow M\rightarrow L\rightarrow 0$ be an exact
sequence of $R$-modules.
\begin{enumerate}
\item
\label{Ass1-exact-sequ}
$\Ass_1(N)\subseteq\Ass_1(M)\subseteq\Ass_1(N)\cup\Ass_1(L)$.
\item
\label{Ass0-exact-sequ}
$\Ass_0(N)\subseteq\Ass_0(M)\subseteq\Ass_0(N)\cup\Ass_0(L)$.
\end{enumerate}
\end{cor}
\begin{pf}
To simplify the notation we may assume that $N\subseteq M$ and
$L=M/N$.\\
Since the annihilator of an element $x\in N$ is the same as the
annihilator of $x$ regarded as an element of $M$ it is trivial that
$\Ass_1(N)\subseteq\Ass_1(M)$ and $\Ass_0(N)\subseteq\Ass_0(M)$.
So all we have to show is that each $\frak p\in\Ass_i(M)$ lies in
$\Ass_i(N)$ or in $\Ass_i(M/N)$ for $i=0,1$.
\\[1ex]
(\ref{Ass1-exact-sequ}): If $\frak p\in\Ass_1(M)$ then by
Proposition \ref{min-prime-ideal-of Ann Rpx}
there is an $x\in M$ such that $R_{\frak p}\cdot x$ is
coprimary for $\frak p\cdot R_{\frak p}$.\\
$1^{st}$ case: $R_{\frak p}\cdot x\cap N_{\frak p}=(0)$. Then $R_{\frak
p}\cdot x\cong R_{\frak p}\cdot x+N_{\frak p}/N_{\frak p}
\subseteq(M/N)_{\frak p}$, which shows that $(M/N)_{\frak p}$ contains a
$\frak p\cdot R_{\frak p}$-coprimary submodule. Then by
Theorem \ref{Ass-M-noetherian}, (\ref{p-coprim-p-Ass})\
$\frak p\cdot R_{\frak p}\in\Ass_1\left((M/N)_{\frak p}\right)$
and then by Proposition \ref{Ass1-up-down}\  $\frak p\in\Ass_1(M/N)$.\\
$2^{nd}$ case: $U:=R_{\frak p}\cdot x\cap N_{\frak p}\ne(0)$. Then by
Proposition \ref{submodule-coprimary-is coprimary}
$U$ is $\frak p\cdot R_{\frak p}$-coprimary as a
submodule $\ne (0)$ of $R_{\frak p}\cdot x$. Since $U\subseteq N_{\frak p}$
it follows from Theorem \ref{Ass-M-noetherian}, (\ref{p-coprim-p-Ass})\
that $\frak p\cdot R_{\frak p}\in\Ass_1(N_{\frak p})$ and then by
Proposition \ref{Ass1-up-down}\  $\frak p\in\Ass_1(N)$.
\\[1ex]
(\ref{Ass0-exact-sequ}):
If $\frak p\in\Ass_0(M)$ there exists a submodule $U=R\cdot x\cong R/\frak p$
of $M$.\\
$1^{st}$ case: $U\cap N=(0)$. Then $U\cong U+N/N\subseteq M/N$ and
therefore $\frak p\in\Ass_0(M/N)$.\\
$2^{nd}$ case: $U\cap N\ne (0)$. Let $0\ne y\in U\cap N$. Since
$U\cong R/\frak p$ and $\frak p$ is a prime ideal, $\Ann_R(y)=\frak p$, and
therefore $\frak p\in\Ass_0(N)$.
\end{pf}
\begin{rem} While $\Ass(N)\subseteq\Ass(M)$ by Corollary
\ref{Ass-submodule}, in general\\
$\Ass(M)\nsubseteq\Ass(N)\cup\Ass(L)$,
even if the exact sequence splits,\\
as is shown in Example \ref{Ass-exact-sequ-false}.
\end{rem}
\section{The Support of a Module}
\begin{defn}
Let $M$ be an $R$-module.
$$
\Supp(M):=\{\frak p\mid\frak p\text{ prime ideal of }M
\text{ with }M_{\frak p}\ne0\}
$$
\end{defn}
We summarize the basic properties of $\Supp$:
\begin{rem}
\label{Supp-basics}
\strut
\begin{enumerate}
\item
\label{M-not-0-Supp-not-empty}
$M\ne0\Longleftrightarrow \Supp(M)\ne\emptyset$
\item
\label{Supp-exact-sequence}
If $0\rightarrow N\rightarrow M\rightarrow L\rightarrow 0$ is an
exact sequence of $R$-modules then
$\Supp(M)=\Supp(N)\cup\Supp(L)$.
\item
\label{Supp-bigger-ideal}
If $\frak p\in\Supp(M)$ and $\frak p'$ is a prime ideal with
$\frak p'\supseteq\frak p$ then $\frak p'\in\Supp(M)$.
\item
\label{Supp-and-Ann}
$\frak p\in\Supp(M)\Longrightarrow\frak p\supseteq\Ann_R(M)$\\
(but the converse is not true in general: There may be prime ideals
containing $\Ann_R(M)$ which do not belong to $\Supp(M)$ as is
shows in Example \ref{Ann-not-in-Supp}.)\\
If $M$ is finitely generated then each $\frak p$ containing
$\Ann_R(M)$ also belongs to $\Supp(M)$.
\item
\label{Supp-of-sums}
$\Supp\sum\limits_{i\in I} N_i=\bigcup\limits_{i\in I}\Supp(N_i)$\\
for arbitrary families of submodules $N_i$ of $M$.
\end{enumerate}
\end{rem}
\begin{pf}
(\ref{M-not-0-Supp-not-empty})
Trivially if $M=(0)$ the $M_{\frak p}=(0)$ for all $\frak p$.\\
Conversely: If $M\ne(0)$ let $0\ne x\in M$. Then $\Ann_R(x)\subsetneq
R$, and so there is a prime ideal $\frak p$ with
$\Ann_R(x)\subseteq\frak p$. Then $0\ne\frac{x}{1}\in M_{\frak p}$ and
therefore $\frak p\in\Supp(M)$.
\\[1ex]
(\ref{Supp-exact-sequence})
For each $\frak p$ the sequence
$0\rightarrow N_{\frak p}\rightarrow M_{\frak p}\rightarrow
L_{\frak p} \rightarrow 0$ is exact. Therefore $M_{\frak p}\ne(0)$
iff $N_{\frak p}\ne(0)$ or $L_{\frak p}\ne(0)$.
\\[1ex]
(\ref{Supp-bigger-ideal})
Because of $M_{\frak p}=\left(M_{\frak p'}\right)_{\frak p\cdot
R_{\frak p'}}$ from $M_{\frak p}\ne(0)$ follows $M_{\frak p'}\ne(0)$.
\\[1ex]
(\ref{Supp-and-Ann})
If $\frak p\nsupseteq\Ann_R(M)$ there is an $s\in\complement\frak p
\cap\Ann_R(M)$, so that $s\cdot M=(0)$ and therefore
$M_{\frak p}=(0)$. For the converse if $M$ is finitely generated see
\ref{Supp-and-Ass}
\\[1ex]
(\ref{Supp-of-sums})
This follows immediately from the fact that
$M_{\frak p}=\sum\limits_{i\in I}{N_i}_{\frak p}$
\end{pf}
There is a close connection between the support of a module and its
associated prime ideals:
\begin{prop}
\label{Supp-and-Ass}
Let $M$ be an $R$-module and $\frak p$ a prime ideal of $R$. Consider
the following conditions:
\begin{enumerate}
\item
\label{p-in-Supp}
$\frak p\in\Supp(M)$
\item
\label{p-Ass-in-p}
$\frak p$ contains a prime ideal of $\Ass(M)$.
\item
\label{p-ess-for submod}
$\frak p$ is essential for a submodule of $M$.
\item
\label{p-ess1-for submod}
$\frak p$ is essential of the first kind for a submodule of $M$.
\item
\label{p-contains-Ann}
$\frak p\supseteq\Ann_R(M)$.
\end{enumerate}
Then {\fontshape{n}\selectfont(\ref{p-in-Supp})--(\ref{p-ess1-for
submod})}  are equivalent,
and {\fontshape{n}\selectfont(\ref{p-contains-Ann})} follows from them.\\
If $M$ is finitely generated then
{\fontshape{n}\selectfont(\ref{p-in-Supp})--(\ref{p-contains-Ann})}
 are equivalent.
\end{prop}
\begin{pf}
(\ref{p-in-Supp}) $\Rightarrow$ (\ref{p-Ass-in-p}):
Let $\frak p\in\Supp(M)$. Then $M_{\frak p}\ne(0)$ and therefore
$\Ass(M_{\frak p})\ne\emptyset$ by Remark \ref{Ass1-rems}
(\ref{M-not-zero-Ass1}). Let $\frak P'\in\Ass(M_{\frak p})$ and
$\frak p':=\frak P'\cap R$. Then $\frak p'\in\Ass(M)$ by
Corollary \ref{Ass-up-down} and $\frak p'\subseteq\frak p$.
\\[1ex]
(\ref{p-Ass-in-p}) $\Rightarrow$ (\ref{p-in-Supp}):
Let $\frak p'\in\Ass(M)$ and $\frak p'\subseteq\frak p$. Then by
Corollary \ref{p-Ass-Mp-not-0}
$M_{\frak p'}\ne(0)$. But then a fortiori $M_{\frak p}\ne(0)$.
\\[1ex]
(\ref{p-in-Supp}) $\Rightarrow$ (\ref{p-ess1-for submod}):
Let $\frak p\in\Supp(M)$. Then $M_{\frak p}\ne(0)$. Let
$0\ne y\in M_{\frak p}$, and
$\widetilde N:=\frak p\cdot R_{\frak p}\cdot y$. Then
$y\notin\widetilde N$ by Krull-Nakayama. Let $N:=\widetilde N\cap
M$, and therefore $\widetilde N=N_{\frak p}$.
We will show that $\frak p\in\Ass_1(M/N)$:\\
Because of Proposition \ref{Ass1-up-down} it is enough to show that
$\frak p\cdot R_{\frak p}\in\Ass_1\left((M/N)_{\frak p}\right)=
\Ass_1(M_{\frak p}/\widetilde N)$:\\
If we denote by $\bar y$ the residue class of $y$ mod $\widetilde N$
we have $\bar y\ne0$, since $y\notin\widetilde N$, but
$\frak p\cdot R_{\frak p}\cdot\bar y=0$ and therefore
$\frak p\cdot R_{\frak p}\subseteq\Ann_{R_{\frak p}}(\bar y)$. But
$\frak p\cdot R_{\frak p}$ is the maximal ideal of $R_{\frak p}$ and
so equality holds. That means that
$\frak p\cdot R_{\frak p}\in\Ass_0(M_{\frak p}/\widetilde N)
\subseteq\Ass_1(M_{\frak p}/\widetilde N)$.
\\[1ex]
(\ref{p-ess1-for submod}) $\Rightarrow$ (\ref{p-ess-for submod}) is
trivial since $\Ass_1(M/N)\subseteq\Ass(M/N)$.
\\[1ex]
(\ref{p-ess-for submod}) $\Rightarrow$ (\ref{p-in-Supp}):
If $\frak p\in\Ass(M/N)$ then by Corollary \ref{p-Ass-Mp-not-0}
$\left(M/N\right)_{\frak p}\ne(0)$. Because of
$\left(M/N\right)_{\frak p}=M_{\frak p}/N_{\frak p}$ we have
a fortiori $M_{\frak p}\ne (0)$.
\\[1ex]
(\ref{p-in-Supp}) $\Rightarrow$ (\ref{p-contains-Ann}) holds by
Remark \ref{Supp-basics} (\ref{Supp-and-Ann}).\\
Now let $M$ be finitely generated. Then we show\\
(\ref{p-contains-Ann}) $\Rightarrow$ (\ref{p-in-Supp}):
Let $\frak p\supseteq\Ann_R(M)$. For a finitely generated $R$-module
one has\\
$\Ann_{R_{\frak p}}(M_{\frak p})=R_{\frak p}\cdot\Ann_R(M)$ (while in
general only $\Ann_{R_{\frak p}}(M_{\frak p})\supseteq R_{\frak
p}\cdot\Ann_R(M)$ (Remark \ref{transporteurs}
(\ref{transp-quotmodule}))).  Therefore
$\Ann_{R_{\frak p}}(M_{\frak p})\subseteq
\frak p\cdot R_{\frak p}\ne R_{\frak p}$ and therefore\\
$M_{\frak p}\ne(0)$.
\end{pf}
 From (\ref{p-Ass-in-p})$\Rightarrow$(\ref{p-in-Supp}) we get:
\begin{cor}
\label{Ass-in-Supp}
$$
\Ass(M)\subseteq\Supp(M)
$$
\end{cor}
\begin{cor}
\label{Supp-R-mod-a}
Let $\frak a$ be an ideal of $R$. Then
$$
\Supp(R/\frak a)=\{\frak p\mid \frak p\supseteq\frak a\}.
$$
\end{cor}
\begin{pf}
$\frak a=\Ann_R(R/\frak a)$ and $M:=R/\frak a$ is a finitely generated
$R$-module. Proposition \ref{Supp-and-Ass}
(\ref{p-in-Supp}) $\Leftrightarrow$ (\ref{p-contains-Ann})
gives the corollary.
\end{pf}
\begin{cor}
\label{min-in-Ass-and-Supp}
\strut
\begin{enumerate}
\item
\label{min-Supp-min-Ass}
$\frak p'$ is minimal in $\Supp(M)\Leftrightarrow\frak p'$ is minimal
in $\Ass(M)$.
\item
\label{M-finite-ex-min-prime}
If $M$ is finitely generated each $\frak p\in\Supp(M)$ contains a
minimal\\
$\frak p'\in\Supp(M)$.\\
If $M$ is not finitely generated there
may be no minimal elements in $\Supp(M)$.
(See Example \ref{no-minimal-primes-in-Supp}.)
\end{enumerate}
\end{cor}
\begin{pf}
(\ref{min-Supp-min-Ass}):
Let $\frak p'$ be minimal in $\Supp(M)$. Then by Proposition
\ref{Supp-and-Ass} (\ref{p-in-Supp}) $\Rightarrow$ (\ref{p-Ass-in-p})
there is a $\frak p''\in\Ass(M)$ with $\frak p'\supseteq\frak p''$,
which by Corollary \ref{Ass-in-Supp}
lies in $\Supp(M)$ and therefore $\frak p'=\frak p''$, because of the
minimality of $\frak p'$.
So we get $\frak p'\in\Ass(M)$ and $\frak p'$ is also minimal in
$\Ass(M)$, because of $\Ass(M)\subseteq\Supp(M)$.\\
Conversely, by the same arguments we see that a minimal prime ideal
of $\Ass(M)$ is also minimal in $\Supp(M)$.
\\[1ex]
(\ref{M-finite-ex-min-prime}):
By Proposition \ref{Supp-and-Ass}\
$\Supp(M)=\{\frak p\mid\frak p\supseteq\Ann_R(M)\}$.
Now each $\frak p\supseteq\Ann_R(M)$ contains a
$\frak p'\supseteq\Ann_R(M)$, which is minimal among the prime ideals
containing $\Ann_R(M)$ and therefore minimal in $\Supp(M)$.
\end{pf}
\section{The Radical of a Submodule}
\begin{defn}
Let $N$ be a proper submodule of $M$. We define the
``radical of $N$ in $M$'' as
$$
\frak r_M(N):=\{r\mid r\in R,\ r \text{ nilpotent for } M/N\}.
$$
\end{defn}
\begin{rem}
\label{basics-radical}
\strut
\begin{enumerate}
\item
$\frak r_M(N)$ is an ideal of $R$.
\item
\label{radN-rad0}
$\frak r_M(N)=\frak r_{M/N}((0))$.
\item
$\frak r_M((0))\supseteq\Ann_R(M)$.
\item
If $\frak a$ is an ideal of $R$ then
$\frak r_R(\frak a)\supseteq\frak a$.
\item
$\frak r_R((0))=\{\text{nilpotent elements of R}\}$ is the
``nil-radical'' of $R$.
\end{enumerate}
\end{rem}
\begin{prop}
\label{nilpotents-intersection-of Ass}
$$
\frak r_M(N)=\bigcap\limits_{\frak p\in\Supp(M/N)}\frak p
            =\bigcap\limits_{\frak p\in\Ass(M/N)}\frak p
            =\bigcap\limits_{\frak p\in\Ass_1(M/N)}\frak p
$$
\end{prop}
\begin{pf}
Since $\Ass_1(M/N)\subseteq\Ass(M/N)$ and
each $\frak p\in\Ass(M/N)$ contains a $\frak p'\in\Ass_1(M/N)$ one has
$\bigcap\limits_{\frak p\in\Ass(M/N)}\frak p
=\bigcap\limits_{\frak p\in\Ass_1(M/N)}\frak p$.
\\[1ex]
Because of Corollary \ref{Ass-in-Supp}
we have $\Ass(M/N)\subseteq\Supp(M/N)$ and each\\
$\frak p\in\Supp(M/N)$ contains a $\frak p'\in\Ass(M/N)$.\\
Therefore $\bigcap\limits_{\frak p\in\Supp(M/N)}\frak p
=\bigcap\limits_{\frak p\in\Ass(M/N)}\frak p$.
\\[1ex]
Because of Remark \ref{basics-radical} (\ref{radN-rad0}) we may now assume
that $N=(0)$. Then $M\ne(0)$.
Let $r\in\frak r_M((0))$. Then for each $x\in M$ there is a
$\nu\in\Bbb N$ with $r^\nu\cdot x=0$ and therefore $M_S=(0)$ for each
multiplicatively closed set $S$ which contains $r$. But for each
$\frak p\in\Supp(M)$ by definition $M_{\complement\frak p}\ne (0)$
and therefore $r\in\frak p$ for all $\frak p\in\Supp(M)$, i.e.
$\frak r_M((0))\subseteq\bigcap\limits_{\frak p\in\Supp(M)}\frak
p$.\\
Conversely: Let $r\in\bigcap\limits_{\frak p\in\Supp(M)}\frak p$ and
let $S:=\{r^\nu\mid\nu=0,1,2,\dots\}$. We will show that $M_S=(0)$,
which means that for each $x\in M$ there is a $\nu\in\Bbb N$ with
$r^\nu\cdot x=0$ and so $r\in\frak r_M((0))$:\\
$1^{st}$ case: $R_S=(0)$. Then also $M_S=(0)$, since we always assume
$M$ to be a unitary $R$-module and in $R=(0)$ the $0$ is the unit
element.\\
$2^{nd}$ case: $R_S\ne(0)$, hence $1\ne0$ in $R_S$. If $M_S\ne(0)$
there would be an $x\in M$ with $\frac{x}{1}\ne0$ in $M_S$, hence
$\Ann_{R_S}\left(\frac{x}{1}\right)\ne R_S$. Then there
would be a prime ideal $\frak P$ of $R_S$ with
$\frak P\supseteq\Ann_{R_S}\left(\frac{x}{1}\right)$.
Let $\frak p:=\frak P\cap R$. Then $\frak p\cap S=\emptyset$, hence
$r\notin\frak p$.\\
But $M_{\frak p}=\left(M_S\right)_{\frak P}\ne(0)$, since
$\frak P\supseteq\Ann_{R_S}\left(\frac{x}{1}\right)$. Therefore
$\frak p\in\Supp(M)$ and so $r\in\frak p$, contradiction!
\end{pf}
\begin{cor}
\label{coprimary-iff-Assp}
$$
\Ass(M)=\{\frak p\}\Longleftrightarrow M\text{ is }\frak p\text{-coprimary}.
$$
\end{cor}
\begin{pf}
If $M$ is $\frak p$-coprimary then $\Ass(M)=\{\frak p\}$ by Remark
\ref{coprimary-Ass}.\\
Conversely: If $\Ass(M)=\{\frak p\}$ then by Theorem
\ref{ass-zero-divisors}
$\frak p$ is the set of all zero divisors for $M$ and by
Proposition \ref{nilpotents-intersection-of Ass}
each element of $\frak p$ is nilpotent for $M$.
\end{pf}
\begin{cor}
$\frak p\in\Supp(M)\Longrightarrow\frak p\supseteq\frak r_M((0))$\\
But in general the converse is not true
(Example \ref{p-rad-not-Supp}).\\
If $M$ is finitely generated, then the converse holds.
\end{cor}
\begin{pf}
``$\Longrightarrow$'' holds because of Proposition
\ref{nilpotents-intersection-of Ass}.\\
``$\Longleftarrow$'':
Let $\frak p\supseteq\frak r_M((0))$ then $\frak p\supseteq\Ann_R(M)$
because of Remark \ref{basics-radical}. If $M$ is finitely generated
then also $\frak p\in\Supp(M)$
by Proposition \ref{Supp-and-Ass}.
\end{pf}
\begin{cor}
Let $N$ be a proper submodule of $M$ and $M/N$ finitely generated or
$\Ass(M/N)$ finite. (E.g. if there exists a primary decomposition of
$N$ in $M$.) Then
$$
\frak r_M(N)=\bigcap\limits_{\vbox{\hsize1.8cm\noindent
\size{8}{8pt}\selectfont
$\frak p$ minimal in\\
$\Ass(M/N)$}}
\frak p
$$
\end{cor}
\begin{pf}
We will show that in both cases each $\frak p\in\Ass(M/N)$ contains a
$\frak p'$ which is minimal in $\Ass(M/N)$. (Therefore one can
restrict the intersection
$\frak r_M(N)=\bigcap\limits_{\frak p\in\Ass(M/N)}\frak p$ to the
minimal elements of $\Ass(M/N)$.):
If $\Ass(M/N)$ is finite this is trivial. If $M/N$ is finitely
generated then by Corollary \ref{min-in-Ass-and-Supp}
(\ref{M-finite-ex-min-prime})
each $\frak p\in\Ass(M/N)$ contains a
minimal element $\frak p'$ of $\Supp(M/N)$, and by
Corollary \ref{min-in-Ass-and-Supp} (\ref{min-Supp-min-Ass}) this is
also minimal in $\Ass(M/N)$.
\end{pf}
\begin{cor}
If $\frak a$ is an ideal of $R$ then
$$
\frak r_R(\frak a)=\bigcap\limits_{\frak p\supseteq\frak a}\frak p=
\bigcap\limits_{\vbox{\hsize1.8cm\noindent
\size{8}{8pt}\selectfont
$\frak p$ minimal \\
containing $\frak a$}}
\frak p
$$
\end{cor}
\begin{pf}
By Corollary \ref{Supp-R-mod-a}
$\Supp(R/\frak a)=\{\frak p\mid\frak p\supseteq\frak a\}$
and therefore by \ref{nilpotents-intersection-of Ass}
$\frak r_R(\frak a)=\bigcap\limits_{\frak p\supseteq\frak a}\frak p$.
Since each $\frak p\supseteq\frak a$ contains a $\frak p'$ which is
minimal among the prime ideals containing $\frak a$ one can restrict
the intersection to the minimal ones among the prime ideals.
\end{pf}
\begin{cor}
If $R$ is ``reduced'' (i.e. $\frak r_R((0))=(0)$ ) then
$$
\{\text{zero divisors of }R\}=\bigcap\limits_{\vbox{\hsize2.1cm\noindent
\size{8}{8pt}\selectfont
$\frak p$ minimal\\
prime ideal of $R$}}\frak p
$$
\end{cor}
\begin{pf}
By Remark \ref{Ass1-rems} (\ref{essential-primes-for-ideal}) the minimal
prime ideals of $R$ belong to\\
$\Ass(R/(0))=\Ass(R)$, and therefore by
Theorem \ref{ass-zero-divisors} all of their elements are zero divisors for
$R$.\\
Conversely: If $r$ is a zero divisor for $R$ there is an $s\in R$,
$s\ne0$ with $r\cdot s=0\in\frak p$ for all minimal prime ideals of
$R$. But since by hypothesis
$\bigcap\limits_{\vbox{\hsize2.1cm\noindent
\shape{it}\size{8}{8pt}\selectfont
$\frak p$ minimal\\
prime ideal of $R$}}\frak p=\frak r_R((0))=(0)$ there is a minimal
prime ideal $\frak p$ of $R$ with $s\notin\frak p$ and therefore
$r\in\frak p$.
\end{pf}
\section{The (Counter-)Examples}
\begin{ex}[\bf\boldmath $(0)$ indecomposable in $M$ but $M$ not
coprimary]
\label{0-indcomp-not-coprim}
Let\\
$R$ be a rank 2 discrete valuation ring in the sense of Krull
\cite{Krull-Bewertung} with valuation $\nu$ and
value group $\Bbb Z\times\Bbb Z$ (lexicographically ordered).\\
$R$ has three prime ideals $\frak P_2\supset\frak P_1\supset (0)$.
Let $\pi_1,\pi_2\in R$ be elements with $\nu(\pi_2)=(0,1)$ and
$\nu(\pi_1)=(1,0)$. Then $\frak P_2=R\cdot\pi_2$ is a principal
ideal, but $\frak P_1$ is not finitely generated. A generating set
for $\frak P_1$ is
$\{\pi_1/\pi_2^i\mid i=0,1,2,\dots\}$.
{\fontshape{n}\selectfont
(Compare also \cite{Berger-Modul-diskret-ganz-Diff}).}\\
Let $M:=R/R\cdot\pi_1$.\\
We claim that
\begin{enumerate}
\item $(0)$ is indecomposable in $M$ (i.e. $R\cdot\pi_1$ is
indecomposable in $R$), but
\item $M$ is {\em not} coprimary.
\end{enumerate}
\end{ex}
\begin{pf}
(1) Let $R\cdot\pi_1=\frak a\cap\frak b$ with ideals $\frak a,\frak b$
of $R$. If $R\pi_1\subsetneq\frak a$ and $R\pi_1\subsetneq\frak b$
then an element of value $<\nu(\pi_1)=(1,0)$ must be contained in
$\frak a$ and in $\frak b$. Among these values
$(1,-1)=\nu(\pi_1/\pi_2)$ is the biggest. Since $R$ is a valuation
ring $\pi_1/\pi_2\in\frak a\cap\frak b=R\pi_1$. This cannot happen
because the values of all elements of $R\pi_1$ are
$\ge\nu(\pi_1)=(1,0)$. It follows that $R\pi_1$ is indecomposable in
$R$.\\
(2) $\pi_1/\pi_2\notin R\pi_1$ but
$\pi_2\cdot(\pi_1/\pi_2)=\pi_1\in R\pi_1$ $\Rightarrow \pi_2$ is a
zero divisor for $M$.\quad But $\pi_2$ is not nilpotent for $M$;
because for all $i\in\Bbb N$ we have $\pi_2^i\cdot 1\notin R\pi_1$,
since $\nu(\pi_2^i)=(0,i)<(1,0)$. Therefore $M$ is not coprimary.
\end{pf}
\begin{ex}[\boldmath$\Ass(M)\supsetneq\Ass_1(M)$]
\label{Ass-not-Ass1}
Let\\
$R:=k[X_1,X_2,\dots]$ polynomial ring in countably many
indeterminates over a field $k$,\\
$\frak p_i:=(X_1,X_2,\dots,X_i)$\\
$\frak p:=(X_1,X_2,X_3,\dots)$,\\
$M:=\bigoplus\limits_{i=1}^\infty R/\frak p_i =
    \bigoplus\limits_{i=1}^\infty R\cdot e_i$
with $e_i:= 1+\frak p_i\in R/\frak p_i$.\\
$\frak p_i$ and $\frak p$ are prime ideals,\\
$\frak p_i=\Ann_R(e_i)$ and therefore
$\frak p_i\in\Ass_0(M)\subseteq\Ass_1(M)$ for $i=1,2,\dots$.\\
$\frak p=\bigcup\limits_{i=1}^\infty\frak p_i$
and therefore $\frak p\in\Ass(M)$ by Proposition
\ref{union-of-Ass1-primes} (\ref{if-union-is-prime}).\\
But $\frak p\notin\Ass_1(M)$.
\end{ex}
\begin{pf}
We have to show that $\frak p$ is not minimal in the set of all prime
ideals containing the annihilator of an element of $M$:\\
Let $0\ne y\in M$ be arbitrary. Then there exists an $n\in\Bbb N$
with $y\in\bigoplus\limits_{i=1}^n R/\frak p_i$.
Let $r\in\Ann_R(y)$ be an arbitrary element of $\Ann_R(y)$.\\
We show that $r\in\frak p_n$:\\
$y=\sum\limits_{i=1}^n\xi_i\cdot e_i$ with $\xi_i\in R$.\\
$r\cdot y=0\Rightarrow r\cdot\xi_i\cdot e_i=0
\Rightarrow r\cdot\xi_i\in\frak p_i$ for $i=1,\dots n$.
But for an $i_0\in\{1,\dots,n\}$ we have
$\xi_{i_0}\notin\frak p_{i_0}$ because else $y=0$.
Then $r\in\frak p_{i_0}\subseteq\frak p_n$.\\
It follows that
$\Ann_R(y)\subseteq\frak p_n\subsetneq\frak p$ and therefore
$\frak p$ is not minimal among the prime ideals containing
$\Ann_R(y)$.
\end{pf}
One can even find a cyclic $R$-module $M$ with
$\Ass(M)\ne\Ass_1(M)$:
\begin{ex}[\boldmath$M$ cyclic and $\Ass(M)\supsetneq\Ass_1(M)$]
Let\\
$R':=k[X_1,Y_1,X_2,Y_2,\dots]$ the polynomial ring in the
countably many independent indeterminates $X_i,Y_i$, $i=1,2,\dots$
over a field $k$,\\
$\frak a':=\left(X_1\cdot Y_1,Y_1^2,X_2\cdot Y_2,Y_2^2,\dots\right)$
ideal in $R'$,\\
$M:=R:=R'/\frak a'=k[x_1,y_1,x_2,y_2,\dots]$, where the $x_i,y_i$ denote
the residue classes of the $X_i,Y_i$ mod $\frak a'$,\\
$\frak p:=(x_1,y_1,x_2,y_2,\dots)\subset R$.\\
Then $\frak p\in\Ass(M)$ but $\frak p\notin\Ass_1(M)$.
\end{ex}
\begin{pf}
(1) The set
$$
A:=\left\{x_1^{\nu_1}\cdot y_1^{\epsilon_1}\cdots x_n^{\nu_n}\cdot
y_n^{\epsilon_n}\ \big|\
n\in\Bbb N,\
\nu_i\in\Bbb N_0,\
\epsilon_i\in\{0,1\},\
\epsilon_i=0 \text{ if }\nu_i>0\right\}
$$
is a basis for $R$ as a $k$-vector space:\\
Obviously any  polynomial of $R'$ can be reduced mod
$\frak a'$ to a linear combination of monomials
$X_1^{\nu_1}\cdot~Y_1^{\epsilon_1}\cdots%
X_n^{\nu_n}\cdot~Y_n^{\epsilon_n}$  with
$n\in\Bbb N,\
\nu_i\in\Bbb N_0,\
\epsilon_i\in\{0,1\},\
\epsilon_i=0 \text{ if }\nu_i>0$ and coefficients in $k$.
Therefore $A$ is a set of generators for $R$ as a $k$-vector space.\\
On the other hand one sees that by definition of $\frak a'$ every
monomial of an element of $\frak a'$ contains an $X_i\cdot Y_i$ or an
$Y_i^2$, while a linear combination of the monomials
$X_1^{\nu_1}\cdot~Y_1^{\epsilon_1}\cdots%
X_n^{\nu_n}\cdot~Y_n^{\epsilon_n}$  with
$n\in\Bbb N,\
\nu_i\in\Bbb N_0,\
\epsilon_i\in\{0,1\},\
\epsilon_i=0 \text{ if }\nu_i>0$ and coefficients in $k$ never
contains these products. Therefore the elements of the set $A$ are
also linearely independent over $k$.
\\[1ex]
(2) No element of $R\setminus\frak p$ is a zero divisor of $R$:\\
Let $F\in R'$ and $F=\sum_{i=0}^n F_i$ its decomposition into
homogeneous polynomials (with respect to the total degree, all
$X_i,Y_i$ having degree $1$).
Then $F$ represents an element of
$R\setminus\frak p$ modulo $\frak a'$ iff $F_0\ne 0$, and
$F\in\frak a'$ iff all $F_i\in\frak a'$ since $\frak a'$ is generated
by monomials (homogeneous elements).\\
Now let $T=\sum_{i=1}^n T_i\in R'$ represent an element
$t\in R\setminus\frak p$ and let $Z=\sum_{i=1}^m Z_i\in R'$ represent an
arbitrary element $z\in R$ with $t\cdot z=0$, i.e.
$T\cdot Z\in\frak a'$. We may assume that $T_0=1$.\\
$T\cdot Z=\sum\limits_\lambda
\left(\sum\limits_{\nu+\mu=\lambda}T_\mu\cdot Z_\nu\right)$
is a homogeneous decomposition. Therefore by the preceding remark
$\sum\limits_{\nu+\mu=\lambda}T_\mu\cdot Z_\nu\in\frak a'$ for all
$\lambda$.
We show by induction on $\lambda$
that all $Z_i\in\frak a'$ and therefore $z=0$:\\
$\lambda=0$:\quad $Z_0\in\frak a'\cap k$=(0).\\
Now let $Z_0,\dots,Z_n\in\frak a'$. Then from
$Z_{n+1}+Z_n\cdot T_1+\dots+Z_0\cdot T^{n+1}\in\frak a'$ we obtain
$Z_{n+1}\in\frak a'$.
\\[1ex]
(3) From (2) we see that the canonical homomorphism
$M\rightarrow M_{\frak p}$ is injective and so
$\frac{y_1}{1}\cdots\frac{y_n}{1}\ne0$ in $M_{\frak p}$ for all
$n=1,2,\dots$, because $y_1\cdots y_n\ne0$ in $M$ for all $n$ as
elements of a $k$-basis of $M$ according to (1).
\\[1ex]
(4) We now proof that $\frak p\in\Ass(M)$ by showing that each
element of $\frak p\cdot R_{\frak p}$ is a zero divisor for
$M_{\frak p}$ (Proposition \ref{zero-divisors-in-Mp}):
Let $p\in\frak p$ be arbitrary. Then $p$ can be written as
$p=\sum\limits_{\nu_i,\mu_i\ge1}a_{\nu_1,\dots,\nu_n,\mu_1,\dots,\mu_n}
\cdot x_1^{\nu_1}\cdots x_n^{\nu_n}\cdot y_1^{\mu_1}\cdots y_n^{\mu_n}$
with $a_{\nu_1,\dots,\nu_n,\mu_1,\dots,\mu_n}\in k$. Since each
summand $\ne0$ contains at least one $x_i$ or $y_i$ with
$i\in\{1,\dots n\}$, and since $x_i\cdot y_i=0$ and $y_i\cdot y_i=0$,
we obtain $p\cdot y_1\cdots y_n=0$ in $M$ and therefore
$\frac{p}{s}\cdot(\frac{y_1}{1}\cdots\frac{y_n}{1})=0$ in $M_{\frak p}$
for each element $s\in R\setminus\frak p$. By (3)
$\frac{y_1}{1}\cdots\frac{y_n}{1}\ne0$ in $M_{\frak p}$, and
therefore $\frac{p}{s}$ is a zero divisor for $M_{\frak p}$.
\\[1ex]
(5) $\frak p\notin\Ass_1(M)$:\\
Proof (indirect): If $\frak p$ was minimal among the prime ideals
containing $\Ann_R(z)$ for a $z\in M$ then by Proposition
\ref{min-prime-ideal-of Ann Rpx} each element of $\frak p\cdot
R_{\frak p}$  would be nilpotent for $R_{\frak p}\cdot z$.
According to (1) $z$ has a representation by the $k$-basis $A$:\\
$z=\sum a_{\nu_1,\dots,\nu_n,\epsilon_1,\dots,\epsilon_n}\cdot
x_1^{\nu_1}\cdot y_1^{\epsilon_1}\cdots x_n^{\nu_n}\cdot
y_n^{\epsilon_n}$ with $\nu_i\in\Bbb N_0$,\ $\epsilon_i\in\{0,1\}$,\
$\epsilon_1=0$ if $\nu_i>0$, and
$a_{\nu_1,\dots,\nu_n,\epsilon_1,\dots,\epsilon_c}\in k$,
not all of them $=0$.\\
Let $m$ be a natural number $m>n$ and $\lambda\in\Bbb N$ arbitrary.\\
Then
$x_m^{\lambda}\cdot z=
\sum a_{\nu_1,\dots,\nu_n,\epsilon_1,\dots,\epsilon_c}\cdot
x_1^{\nu_1}\cdot y_1^{\epsilon_1}\cdots x_n^{\nu_n}\cdot
y_n^{\epsilon_n}\cdot x_m^\lambda$ is again a representation of
$x_m^\lambda\cdot z$ by the basis $A$ and therefore
$x_m^\lambda\cdot z\ne0$,
since not all of the coefficients are $0$.
By (3) it follows that
$(\frac{x_m}{1})^\lambda\cdot\frac{z}{1}\ne0$ in $M_{\frak p}$.
But $x_m\in\frak p$ and so $\frac{x_m}{1}$ is an element of
$\frak p\cdot R_{\frak p}$ which is not nilpotent for
$R_{\frak p}\cdot z$.
\end{pf}
\begin{ex}[{\boldmath$\Ass(M)=\Ass_1(M)\supsetneq\Ass_0(M)
=\emptyset$}]%
\label{Ass=Ass1-not-Ass0}
$M$ cyclic and $(0)$ has a primary decomposition in $M$: Let\\
$R$ a valuation ring with value group $\Gamma=\Bbb Q$ or
$\Gamma=\Bbb R$ (rank one, non discrete),\\
$\nu$ the (additive) valuation of $R$,\\
$\frak P=\{z\mid z\in R,\ \nu(z)>0\}$ the maximal ideal of $R$,\\
$\frak a=\{z\mid z\in R,\ \nu(z)\ge 1\}$,\\
$M:=R/\frak a$.\\We show that:
\begin{enumerate}
\item
\label{zero-primary-dec}
$(0)$ is $\frak P$-primary in $M$. (Therefore $(0)$ has as primary
decomposition in $M$.)
\item
\label{noAnn}
There is no $x\in M$ with $\frak P=\Ann_R(x)$.
\end{enumerate}
It follows from {\shape{n}\selectfont(\ref{zero-primary-dec})} by
Remark \ref{coprimary-Ass} that $\Ass(M)=\{\frak P\}$ and from
{\shape{n}\selectfont(\ref{noAnn})} by definition that $\frak P\notin\Ass_0(M)$
and therefore
$\Ass_0(M)=\emptyset$ since $\Ass_0(M)\subseteq\Ass(M)$.
\end{ex}
\begin{pf} (\ref{zero-primary-dec}): Let $p\in\frak P$ and
$a\in\frak a$. Then there is an $n\in\Bbb N$ with
$n\cdot\nu(p)\ge\nu(a)$, hence $p^n\in R\cdot a\subseteq\frak a$
and therefore $p^n\cdot M=(0)$, i.e. each element of $\frak P$ is
nilpotent for $M$. But since $\frak P$ is the maximal ideal of $R$
all zero divisors for $M$ lie in $\frak P$ and so
it follows that $(0)$ is $\frak P$-primary in $M$.
\\[1ex]
(\ref{noAnn}): If there was an $x\in M$ with $\Ann_R(x)=\frak P$
there would be a representative $z\in R\setminus\frak a$ with
$p\cdot z\in\frak a$ for all $p\in\frak P$, i.e.
$\nu(z)+\nu(p)\ge1$ for all $p\in\frak P$.
Now $\Gamma\supseteq\Bbb Q$, and therefore for each $n\in\Bbb N$
there is a $p_n\in\frak P$ with $\nu(p_n)=\frac{1}{n}$. Then
$\nu(z)+\frac{1}{n}\ge1$ for all $n\in\Bbb N$ and hence
$\nu(z)\ge 1$, i.e. $z\in\frak a$ against our assumption.
\end{pf}
\begin{ex}[{\size{11}{12pt}\selectfont\boldmath $R$ noetherian (local),
$\Ass(M)\supsetneq\Ass_1(M)$($=\Ass_0(M)$)}]
\label{R-noeth-Ass-not-Ass1}\strut\\
Let $R:=k[X,Y]_{(X,Y)}$ localization of the polynomial ring in $X$ and
$Y$ over a field $k$,\\
$\cal P:=\{R\cdot p\mid p\in R,\ R\cdot p\text{ prime ideal of }R\}$,\\
$M:=\bigoplus\limits_{R\cdot p\in\cal P}R/R\cdot p$.\\
Then\\
$\Ass_1(M)=\{R\cdot p\mid R\cdot p\in\cal P\}$,\\
$\frak m:=(X,Y)\in\Ass(M)\setminus\Ass_1(M)$
\end{ex}
\begin{pf}
Denote $U_p:=R/R\cdot p = R\cdot e_p$. Then $\Ann_R(e_p)=R\cdot p$,
hence $R\cdot p\in\Ass_0(M)\subseteq\Ass_1(M)$.\\
Conversely let
$0\ne\xi\in M$ arbitrary and $z\in\Ann_R(\xi)$.
$\xi=\sum\xi_p\cdot e_p$ with
$\xi_{p_0}\cdot e_{p_0}\ne0$ for some $R\cdot p_0\in\cal P$, i.e.
$\xi_{p_0}\notin R\cdot p$.
Now from $z\cdot \xi=0$ we get $z\cdot \xi_{p_0}\cdot e_{p_0}=0$, i.e.
$z\cdot\xi_{p_0}\in R\cdot p_0$ and therefore $z\in R\cdot p_0$
since $\xi_{p_0}\notin R\cdot p_0$ and $R\cdot p_0$ is a prime ideal.
It follows that $\Ann_R(\xi)\subseteq R\cdot p_0$. Since the only non
principal prime ideal of $R$ is $\frak m$, which contains all the
$R\cdot p$, it follows that the minimal
elements among the prime ideals containing the annihilator of an
element of $M$ are the  principal prime ideals
$R\cdot p$. So we obtain that
$\Ass_1(M)=\{R\cdot p\mid R\cdot p\in\cal P\}$.\\
Then $\frak m\notin\Ass_1(M)$. But $\frak m\in\Ass(M)$:\\
In view of Proposition \ref{zero-divisors-in-Mp} we only must show
that each element of
$\frak m$ is a zero divisor for $M$ ($R_{\frak m}=R$ !). To show this
let $z\in\frak m$ be arbitrary. $R$ being a UFD there is a prime
element $p$ and a $z_1\in R$ with $z=z_1\cdot p$. But then
$z\cdot e_p=0$ and therefore $z$ is a zero divisor for $M$.
(Another way of showing $\frak m\in\Ass(M)$ would be to use
Proposition \ref{union-of-Ass1-primes} (\ref{if-union-is-prime}).)
\end{pf}
\begin{ex}[\boldmath$\Ass(N\oplus L)\nsubseteq\Ass(N)\cup\Ass(L)$]
\label{Ass-exact-sequ-false}\strut\\
Let $R$,\  $\frak m$, and $M$ be the same as in Example
\ref{R-noeth-Ass-not-Ass1},\\
$N:=R/R\cdot X$,\\
$L:=\bigoplus\limits_{\size{8}{6pt}\selectfont
                     \begin{array}{c}
                     R\cdot p\in\cal P\\
                    R\cdot  p\ne R\cdot X
                      \end{array}}%
\hbox to 6mm{\hss$R/R\cdot p$}$.
\\[1ex]
Then we have by definition\\
$M=N\oplus L$ and therefore we have the splitting exact sequence\\
$0\rightarrow N\rightarrow M\rightarrow L\rightarrow 0$
\\[1ex]
We show
\begin{enumerate}
\item
\label{m-in-AssM}
$\frak m\in\Ass(M)$
\item
\label{m-notin-AssN}
$\frak m\notin\Ass(N)$
\item
\label{m-notin-AssL}
$\frak m\notin\Ass(L)$
\end{enumerate}
\end{ex}
\begin{pf}
(\ref{m-in-AssM}):
This was already shown in Example \ref{R-noeth-Ass-not-Ass1}.
\\[1ex]
(\ref{m-notin-AssN}):
Since $X$ is a prime element of $R$ and $N=R/R\cdot X$ we see that
$N$ is $R\cdot X$-coprimary and so $\Ass(M)=\{R\cdot X\}\not\ni\frak m$.
\\[1ex]
(\ref{m-notin-AssL}):
By definition of $L$ it is obvious that $X$ is not a zero divisor for
$L$ and therefore $\frak m\notin\Ass(L)$, because $X\in\frak m$.
\end{pf}
\begin{ex}
[{\bf\boldmath$\frak p\supseteq\frak r_M((0))$ but
$\frak p\notin\Supp(M)$}]
\label{p-rad-not-Supp}Let\\
$R:=\Bbb Z$,\\
$M:=\bigoplus
\limits_{\size{8}{9pt}\selectfont
          \begin{array}{c}
                0\ne (p)\\
               \text{ prime ideal}
          \end{array}}%
\hbox to.3cm{\hss$\Bbb Z/(p)$}$.\quad Then\\
$\Supp(M)=\{(p)\mid 0\ne p\text{ prime element in }\Bbb Z\}$,
but\\
$\frak r_M((0))=(0)$.
\\[1ex]
So $\frak p:=(0)\supseteq\frak r_M((0))$, but
$\frak p\notin\Supp(M)$.
\end{ex}
\begin{pf}
$(0)\notin\Supp(M)$, because $M$ is a torsion module and
$R_{(0)}=\Bbb Q$ is a field. Therefore $M_{(0)}=(0)$.\\
$\frak r_M((0))=(0)$: For each $0\ne n\in\Bbb Z$ there is a prime
element $p$ with $p\nmid n$ and therefore $p\nmid n^\nu$ for all
$\nu\in\Bbb N$. So $n^\nu\cdot e_p\ne0$ for all $\nu\in\Bbb N$, with
$e_p:=1+(p)$ in the summand $\Bbb Z/(p)$ of $M$.
Therefore $n\notin\frak r_M((0))$, i.e. $\frak r_M((0))=(0)$.
\end{pf}
\begin{ex}[\bf\boldmath$\frak p\supseteq\Ann_R(M)$, but
$\frak p\notin\Supp(M)$]
\label{Ann-not-in-Supp}\strut\\
Let:
$R:=\Bbb Z$,\quad
$M:=\Bbb Q/\Bbb Z$.\\
Then $\Ann_R(M)=(0)$, but
$(0)\notin\Supp(M)$
since $R_{(0)}=\Bbb Q$ is a field and $M$ is a torsion module.
\end{ex}
\begin{ex}
[{\bf\boldmath No minimal elements in $\Supp(M)$}]
\label{no-minimal-primes-in-Supp}
Let\\
$R:=k[X_1,X_2,\dots]$ a polynomial ring in countably many
indeterminates over a field $k$,\\
$\frak p_i:=(X_i,X_{i+1},\dots)$,\\
$M:=\bigoplus\limits_{i=1}^\infty R/\frak p_i$.
\\[1ex]
There are no minimal elements in $\Supp(M)$
\end{ex}
\begin{pf}
Let $M_i:=R/\frak p_i$. By Corollary \ref{Supp-R-mod-a} we have
$\Supp(M_i)=\{\frak p\mid\frak p\supseteq\frak p_i\}$.
Further by Remark \ref{Supp-basics} (\ref{Supp-of-sums})
$\Supp(M)=\bigcup\limits_{i=1}^\infty\Supp(M_i)$. For every
$\frak p\in\Supp(M)$ there is an $i_0$ with
$\frak p\in\Supp(M_{i_0})$ and therefore
$\frak p\supseteq\frak p_{i_0}\supsetneq\frak p_{i_0+1}
\supsetneq\dots$ and all the $\frak p_i\in\Supp(M)$.
So obviously $\frak p$ is not minimal in $\Supp(M)$.
\end{pf}
\begin{ex}[{\size{11}{12pt}\selectfont\bf\boldmath
Essential prime ideals for $\Ann_R(M)$ not in $\Ass(M)$}]
\label{0-primary-decomp-Ann-not ass}\strut\\
Let
$R$ a rank one discrete valuation ring,\quad
$\frak p$ the maximal ideal of $R$,\\
$M:=\Quot(R)/R$.\\
$\Ann_R(M)=(0)$, a prime ideal of $R$.\\
Then\\
$\Ann_R(M)$ has a primary decomposition and $(0)$ is essential for
$\Ann_R(M)$,\\
but\\
$M$ is $\frak p$-coprimary, because each element of $\frak p$ is
nilpotent for $M$ while the elements of $R\setminus\frak p$ are units
of $R$. Then $\Ass(M)=\{\frak p\}\not\ni(0)$
(Remark \ref{coprimary-Ass}).
\end{ex}
\vspace{1cm}
\noindent Author's address:
\\[1ex]
Robert W. Berger\\
Fachbereich Mathematik\\
Universit\"at des Saarlandes\\
Postfach 15 11 50\\
D-66041 Saarbr\"ucken\\
Germany
\\[1ex]
E-mail address:
%\\[1ex]
rberger@@math.uni-sb.de
\end{document}